\documentclass[twocolumn]{autart}
\usepackage{graphicx}
\usepackage{amsmath}
\usepackage{amsfonts}

\begin{document}

\begin{frontmatter}

    \title{Minimum-time lateral interception of a moving target by a Dubins car\thanksref{footnoteinfo}}
    
    \thanks[footnoteinfo]{This paper was not presented at any IFAC meeting.}
    
    \author[Moscow]{Maksim Buzikov}\ead{me.buzikov@physics.msu.ru},
    \author[Moscow]{Andrey Galyaev}\ead{galaev@ipu.ru}
    
    \address[Moscow]{V. A. Trapeznikov Institute of Control Sciences of Russian Academy of Sciences, Moscow, Russia}
              
    \begin{keyword}
        Dubins vehicle; Moving target; Minimum-time interception; Analytic Design; Application of Nonlinear Analysis and Design.
    \end{keyword}
    
    \begin{abstract}
        This paper presents the problem of lateral interception by a Dubins car of a target that moves along an a priori known trajectory. This trajectory is given by two coordinates of a planar location and one angle of a heading orientation, every one of them is a continuous function of time. The optimal trajectory planning problem of constructing minimum-time trajectories for a Dubins car in the presence of a priory known time-dependent wind vector field is a special case of the presented problem. Using the properties of the three-dimensional reachable set of a Dubins car, it is proved that the optimal interception point belongs to a part of an analytically described surface in the three-dimensional space. The analytical description of the surface makes it possible to obtain 10 algebraic equations for calculating parameters of the optimal control that implements the minimum-time lateral interception. These equations are generally transcendental and can be simplified for particular cases of target motion (e.g. resting target, straight-line uniform target motion). Finally, some particular cases of the optimal lateral interception validate developments of the paper and highlight the necessity to consider each of 10 algebraic equations in general case.
    \end{abstract}

\end{frontmatter}

\section{Introduction}

Missions for autonomous ground, aerial and marine vehicles are extremely diverse. In most scenarios, the mission involves a mobile vehicle path/trajectory planning (PP/TP) problem to be solved. Formalization of the problem requires choosing a model of system dynamics that satisfies several requirements at once. First of all, reference trajectories given by model dynamics should be close to feasible ones in practice. Secondly, the computing of these trajectories should be a time-rapid and memory-efficient process for an onboard computer. The analytical results obtained from such models tend to speed up computing of trajectories and reduce memory usage.

Dubins car is a popular model of a non-holonomic vehicle that moves only forward on the plane along trajectories of bounded curvature. The optimal PP-problem of finding a minimum-length path with bounded curvature was for the first time considered by Markov in context of railway track construction \cite{markov1889some}. In \cite{dubins1957curves}, Dubins solved a similar problem, when the orientation of the line is fixed at both ends, in contrast to Markov's statement, where it is arbitrary at the final point and given at the starting one. In \cite{isaacs1965differential} Isaacs considered pursuit-evasion games, and in one of the problems he used equations of dynamics that correspond to the motion with a constant velocity and bounded maneuverability. Isaacs probably was the first who introduced the term "car" for such a model. In modern literature, such a model is called "Dubins car". A complete analytical solution of the Dubins problem is obtained in \cite{pecsvaradi1972optimal,bui1994shortest}. In the case when the distance between initial and terminal points is relatively large it is possible to construct a logical classification scheme for time-rapid computation of the optimal control \cite{shkel2001classification}. A large number of modifications to the Dubins problem has been considered in literature. For example, recent studies devoted to the Dubins PP-problem deal with the case when the trajectory necessarily passes through a fixed intermediate point \cite{chen2019shortest}. Also recently the problems of constructing the shortest path to a circle \cite{chen2020dubins} and the shortest path via a circular boundary \cite{jha2020shortest} for Dubins car have been solved. These analytical results make it affordable to implement efficient onboard algorithms.

The Dubins' model study is not limited to the PP/TP problems. Some theoretically important works are devoted to the study of the Dubins car reachable set. A detailed analysis of the planar evolution of the Dubins car reachable set is given in \cite{cockayne1975plane}. The evolution of a three-dimensional reachable set was studied in \cite{patsko2003three}, and in \cite{patsko2020analiticheskoe}, where a partial analytical description of the evolution of boundary set was obtained. From a practical point of view the reachable sets can be used, for example, to recover the vehicle trajectory from inaccurate measurements \cite{bedin2010restoration} and also in the collision avoidance problem \cite{sahawneh2015airborne}. An analytical description of the planar reachable set made it possible to solve the problem of intercepting a moving target by a Dubins car at an unspecified interception angle \cite{buzikov2021time,zheng2021time}, as well as to clarify necessary conditions for interception along a geodetic line (relaxed Dubins path) \cite{meyer2015dubins}.

This paper solves the problem of minimum-time lateral interception of a moving target by Dubins car when the planar coordinates and the course angles of the target and Dubins car must coincide at the interception moment. The case of unlimited maneuverability, where the criterion of optimality is the sum of the time and the energy spent on maneuvers, is analyzed in \cite{glizer1996optimal}. A large number of works deal with the problem of achieving the final state in the presence of wind instead of intercepting a moving target problem. The case of the presence of constant wind and the case of the straight-line uniform motion of the target are equivalent in the sense that there is a coordinate transformation in which one case passes into another \cite{mcgee2005optimal}. Some properties of the constant wind case solution are investigated in \cite{mcgee2007optimal}. The construction of a numerical scheme for computing the case of time-variable wind can be found in \cite{mcneely2007tour}. The first analytical results for the constant wind case were obtained in the assumption of a sufficient relative distance between the target and Dubins car when considering a selection of aerobiological samples mission \cite{techy2008path}. Later, the same authors obtained analytical results for the general case of the initial target location \cite{techy2009minimum}. The main result of these papers is that for some configurations there are explicit expressions for calculating the parameters of optimal control, and for other scenarios, in the worst case, it is necessary to solve a system of two transcendental algebraic equations with two unknowns. A numerical method for solving a similar system of transcendental equations has been considered in \cite{ayhan2018time}. It is shown in \cite{bakolas2010time} that the system can be simplified to a transcendental equation with one unknown. Analytical results of the problem of achieving by Dubins car the desired state in the presence of constant wind are summarized in \cite{bakolas2013optimal}. In recent paper \cite{mittal2020rapid}, it was shown that by losing optimality, admissible control can be obtained without solving transcendental equations for the constant wind case.

Direct addressing the problem of lateral interception of a target with straight-line uniform motion can be found in \cite{gopalan2016time}. This article deals only with the case of intercepting a target at a right angle under the assumption of a large enough initial distance between the target and Dubins car. The paper \cite{gopalan2017generalized} generalizes this problem to the case of an arbitrary given angle of interception. The paper \cite{manyam2019optimal} provides a solution to a similar problem, but for the case of a target moving in a circle. It is assumed that the intercepting velocity vectors should become collinear. A more complex case of a target moving along racetrack path is discussed in \cite{manyam2020intercepting}.

The statement of the problem considered in this paper generalizes all of the above statements. No restrictions are imposed on the target movement, except for the continuity of the trajectory on time. The area of the practical application of such developments is significantly wide. It includes the tasks of constructing reference trajectories for various controllable moving plants in the presence of wind \cite{schopferer2015performance,coombes2017landing}, the task of refueling an aircraft and other vehicles \cite{burns2007simulation,manyam2019optimal}, the task of landing UAVs on a moving platform \cite{alijani2020autonomous}. The same task can be solved, including wind presence, in military applications \cite{gopalan2016time}, such as tasks of interception and destruction of targets, moving along program trajectories.

This paper is organized as follows: Section \ref{sec:statement} provides the mathematical statement of the problem of lateral interception of a moving target by a Dubins car and the implied terminology. Next Section \ref{sec:optimal_terminal_point} describes the position of the optimal interception point and refines some properties of the three-dimensional reachable set. Then Section \ref{sec:optimal_control} is devoted to obtaining equations for calculating the optimal interception time moment and determining all parameters of the optimal control by the time moment. Section \ref{sec:examples} addresses considering some examples of the optimal interception. Finally, Section \ref{sec:conclusion} gives conclusions and future work.

\section{Problem formulation}\label{sec:statement}

In this section, we formulate the problem of lateral interception of a moving target as an optimal control problem. The subscript $E$ used throughout this paper refers to the target configuration (endpoint for Dubins car). We use the notation $P \in \mathbb{R}^2 \times \mathbb{S}$ for a configuration, which is a planar position and orientation triplet. The space $\mathbb{S}$ consists of the real numbers, but the equality on $\mathbb{S}$ differs from the equality on $\mathbb{R}$. For every $\varphi,\psi \in \mathbb{S}$ the equality $\varphi = \psi$ holds on $\mathbb{S}$ when there exists a number $k \in \mathbb{Z}$ such that the equality $\varphi = \psi + 2\pi k$ holds on $\mathbb{R}$. We also suppose that the absolute value of $\varphi \in \mathbb{S}$ is calculated\footnote{Modulo operation is supposed to be a remainder of division: $\mathrm{mod}(a, b) = a - b\lfloor a/b \rfloor$.} according to the next rule
\begin{equation*}
    |\varphi| =  \min(\mathrm{mod}(\varphi, 2\pi), 2\pi - \mathrm{mod}(\varphi, 2\pi)).
\end{equation*}
We consider throughout this paper the following metric on the space $\mathbb{R}^2 \times \mathbb{S}$
\begin{equation*}
    \rho(P, P_0) = \sqrt{(x - x_0)^2 + (y - y_0)^2 + |\varphi - \varphi_0|^2},
\end{equation*}
where $P = (x, y, \varphi) \in \mathbb{R}^2 \times \mathbb{S}$, $P_0 = (x_0, y_0, \varphi_0) \in \mathbb{R}^2 \times \mathbb{S}$.

\subsection{Dubins' model}
\begin{figure}
    \begin{center}
        \includegraphics[height=7.3cm]{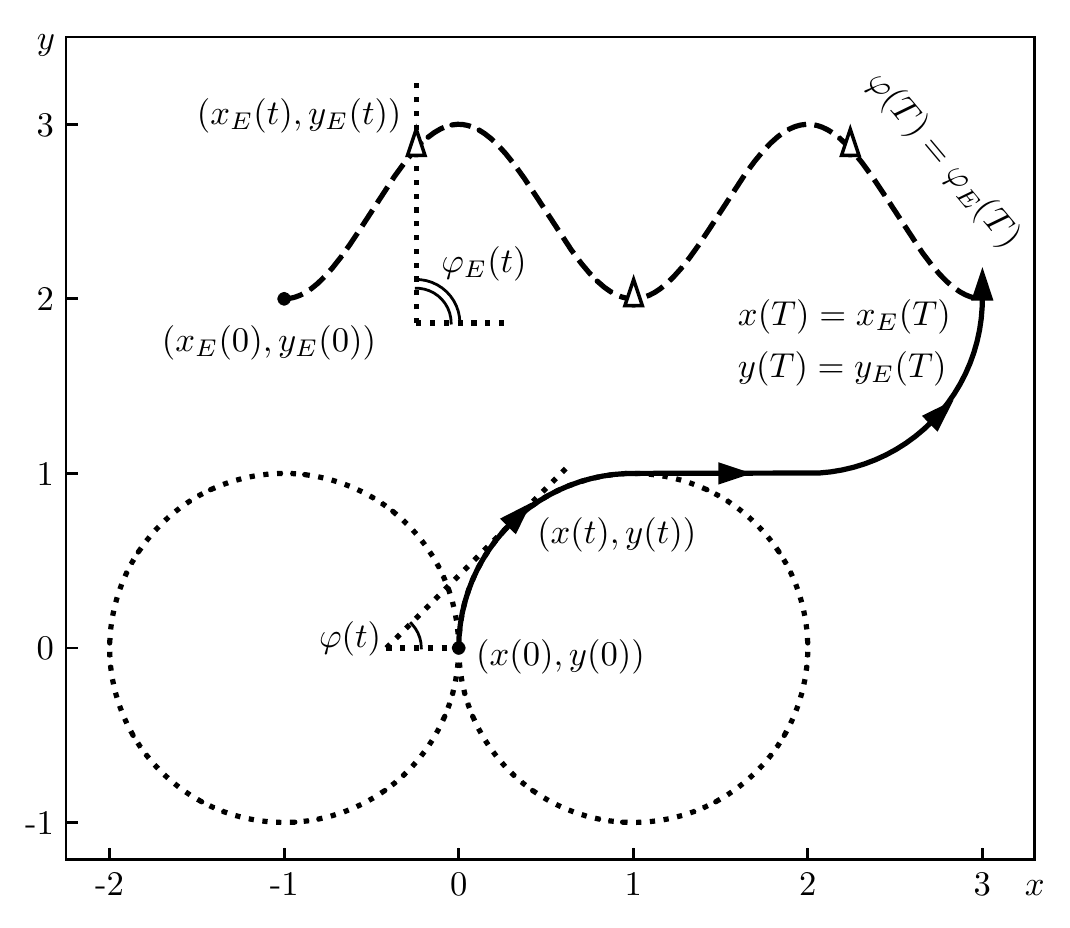}
        \caption{The Dubins car planar path is a solid line. It starts from the configuration $(0, 0, \pi/2)$ and ends on $(3, 2, \pi/2)$. The target planar path is dashed and its orientation is supposed to be constantly upward directed ($\varphi_E(t) = \pi/2$). Dotted circles have a unit radius and the Dubins car is not able to be inward these circles for the first time.}
        \label{fig:problem_statement}
    \end{center}
\end{figure}

Dubins' model describes a vehicle planar motion with constant-speed and limited maneuverability \cite{bui1994shortest}. We suppose that time and length scales are chosen in such a way that the velocity and the minimum turn-radius of the Dubins car are units. That is, the Dubins car motion is described by the system of equations
\begin{equation}\label{eq:dynamics}
    \left\{
    \begin{aligned}
    	&\dot{x} = \cos\varphi,\\
    	&\dot{y} = \sin\varphi,\\
    	&\dot{\varphi} = u,
    \end{aligned}
    \right.
\end{equation}
where $(x(t), y(t)) \in \mathbb{R}^2$ are car Cartesian coordinates (see Fig.~\ref{fig:problem_statement}), $\varphi(t) \in \mathbb{S}$ is the heading of forward velocity (we use counterclockwise system and suppose that x-axis is a reference direction), $u(t) \in [-1, 1]$ is the control input at the time $t$. Without loss of generality, we assume that every considered scenario starts at the time moment $t = 0$ with the initial configuration of the Dubins car being 
\begin{equation}\label{eq:initial_conditions}
    (x(0), y(0), \varphi(0)) = (0, 0, \pi/2).
\end{equation}

\subsection{Minimum-time lateral interception}
The configuration of the target is declared by a continuous vector-function $E = (x_E, y_E, \varphi_E)$, where $(x_E(t), y_E(t)) \in \mathbb{R}^2$ are the target Cartesian coordinates (see Fig.~\ref{fig:problem_statement}), $\varphi_E(t) \in \mathbb{S}$ is the target orientation. Generally speaking, the heading of the target velocity should not be the same as the target orientation. Moreover, the target velocity value may not exist when $(x_E, y_E)$ is not differentiable. Such description of target motion allows considering exotic cases of movement. For example, a resting and rotating point can be described by $E(t) = (x_0, y_0, \alpha t)$.

The lateral interception occurs when car and target configurations coincide with each other
\begin{equation}\label{eq:terminal_conditions}
	(x(T), y(T), \varphi(T)) = (x_E(T), y_E(T), \varphi_E(T)),
\end{equation}
where $T \in \mathbb{R}^+_0$ is an interception time.

We examine the optimal-control problem, where the goal of the Dubins car control is to minimize the interception time $T$, as it is pointed below 
\begin{equation*}
	J[u] \overset{\mathrm{def}}{=} \int\limits^T_0 dt \to \min_{u}.
\end{equation*}
According to Theorem~1 in \cite{cockayne1975plane}, the piecewise-constant control functions are enough to produce interception trajectories in the considered problem. The conditions on the target trajectory for a successful interception will be established further.

\subsection{Background}
Throughout this paper, we use the common notation and some known results relating to the Dubins' model. Similar to \cite{bui1994shortest} we suppose that C denotes an arc of a unit circle and S denotes a straight line segment. Sequences CSC, CCC are assigned to the trajectory consisting only of C, S consecutive parts. For definiteness, the controls raising these trajectories are supposed to be right-continuous piecewise-constant functions \cite{patsko2003three}. If C corresponds to a clockwise or counterclockwise turn, it will be replaced by R or L respectively. So L corresponds to the control $u(t) = 1$ and R corresponds to $u(t) = -1$, but S matches the control action $u(t) = 0$.

Further analysis significantly uses the reachable set properties. Let $\mathcal{R}(t)$ be the three-dimensional reachable set of the Dubins car at time $t \in \mathbb{R}^+_0$. It means that for any point $P \in \mathcal{R}(t)$ there exists an admissible control, the usage of which leads to the point $P$ at time $t$. We note that $\mathcal{R}: \mathbb{R}^+_0 \to \mathbb{R}^2 \times \mathbb{S}$ is a multi-valued mapping. According to general results of optimal control theory \cite{lee1986foundations} the set $\mathcal{R}(t)$ is closed and bounded. $\mathcal{B}(t)$ denotes a boundary of the set $\mathcal{R}(t)$. The main result of \cite{patsko2003three} is that CSC- and CCC- trajectories are enough to reach any point of $\mathcal{B}(t)$ at time $t$, but not all trajectories of these types lead to $\mathcal{B}(t)$, some of them lead to the internal points of $\mathcal{R}(t)$. Let $\mathcal{E}(t)$ be a set of all possible terminations of CSC- and CCC- trajectories at time $t$. It follows from above that
\begin{equation*}
    \mathcal{B}(t) \subset \mathcal{E}(t) \subset \mathcal{R}(t)
\end{equation*}
at any time moment $t \in \mathbb{R}^+_0$. Further, we use the partition $\mathcal{E}(t) = \mathcal{E}_{CSC}(t) \cup \mathcal{E}_{CCC}(t)$, where $\mathcal{E}_{CSC}(t)$ corresponds to the CSC-trajectories terminal points and $\mathcal{E}_{CCC}(t)$ consists of the CCC-trajectories terminal points.

\section{Optimal interception point}\label{sec:optimal_terminal_point}
This section addresses new properties of mappings $\mathcal{B}$, $\mathcal{E}$, $\mathcal{R}$ and contains the analytical description of the mapping $\mathcal{E}$. These properties allow to determine the optimal interception point configuration and to reduce the space of optimal trajectory candidates.

\subsection{Reduction of CSC-, CCC- trajectory classes} 
First of all, we have to succinctly describe controls that allow reaching points of $\mathcal{E}(t)$. Let $\mathbb{B} = \{-1, 1\}$ be a binary set.
\begin{lem}\label{lem:reduced_CSC_CCC}
    Let $s, \sigma \in \mathbb{B}$, $\tau_1 \in [0, 2\pi)$, $\tau_2 \in [\tau_1, +\infty)$. Then any point of $\mathcal{E}_{CSC}(t)$ is attainable using the control
    \begin{equation}\label{eq:implicit_control_CSC}
    	u^{s, \sigma}_{CSC}(t; \tau_1, \tau_2) =
    	\begin{cases}
    		s, & \quad t \in [0, \tau_1),\\
    		0, & \quad t \in [\tau_1, \tau_2),\\
    		\sigma, & \quad t \in[\tau_2, +\infty),
    	\end{cases}
    \end{equation}
    and any point of $\mathcal{E}_{CCC}(t)$ is attainable using the control
    \begin{equation}\label{eq:implicit_control_CCC}
    	u^{s}_{CCC}(t; \tau_1, \tau_2) =
    	\begin{cases}
    		s,& \quad t \in [0, \tau_1),\\
    		-s,& \quad t \in [\tau_1, \tau_2),\\
    		s,& \quad t \in [\tau_2, +\infty),
    	\end{cases}
    \end{equation}
    where $\tau_2 - \tau_1 \in [0, 2\pi)$ for the CCC-trajectory case.
\end{lem}
The times $\tau_1$, $\tau_2$ are called the first and the second switch times. The binary values $s$, $\sigma$ characterize the turning directions. For example, values $s = 1$, $\sigma = -1$ for the control \eqref{eq:implicit_control_CSC} correspond to the LSR-trajectory and $s = -1$ for the control \eqref{eq:implicit_control_CCC} corresponds to the RLR-trajectory. Since the car has a unit velocity the values $\tau_1$, $\tau_2 - \tau_1$ equal to the lengths of the first and second parts of a path.
\begin{pf}
    Assume initially that a CSC-trajectory terminating at some point $P \in \mathcal{E}_{CSC}(t)$ has cycles at the first circle arc. Consider a new trajectory obtained from the previous trajectory by shifting these cycles to the last circle arc. The new trajectory corresponds to the control \eqref{eq:implicit_control_CSC} and terminates at the point $P$. The same reasoning is applicable when we consider the case of CCC-trajectories with possible cycles. These cycles have to be shifted from the first and second circle arcs to the last one. Such a new trajectory corresponds to the control \eqref{eq:implicit_control_CCC} and terminates at the same point as the previous one does.
\end{pf}

Using this lemma, now we may obtain the analytical description of the set $\mathcal{E}(t)$. According to this lemma, for any point $P^{s, \sigma}_{CSC}(t; \tau_1, \tau_2)\in\mathcal{E}_{CSC}(t)$ there exists a control \eqref{eq:implicit_control_CSC} such that corresponding trajectory terminates at this point at time $t$. Now let's integrate equations \eqref{eq:dynamics} using the initial conditions \eqref{eq:initial_conditions} and the control \eqref{eq:implicit_control_CSC}. It gives the following values
\begin{equation}\label{eq:CSC_integrated}
    \begin{split}
        x_{CSC}^{s, \sigma}(t; \tau_1, &\tau_2) = s(\cos\tau_1 - 1 - (\tau_2 - \tau_1)\sin\tau_1)\\
        &+\sigma(\cos(s\tau_1 + \sigma(t - \tau_2)) - \cos\tau_1),\\
        y_{CSC}^{s, \sigma}(t; \tau_1, &\tau_2) = \sin\tau_1 + (\tau_2 - \tau_1) \cos\tau_1\\
        &-\sigma(s\sin\tau_1 - \sin(s\tau_1 + \sigma(t - \tau_2))),\\
        \varphi_{CSC}^{s, \sigma}(t; \tau_1, &\tau_2) = \pi/2 + s\tau_1 + \sigma(t-\tau_2),
    \end{split}
\end{equation}
where we consider that
\begin{equation*}
    \begin{split}
        P^{s, \sigma}_{CSC}(t; \tau_1, \tau_2) \overset{\mathrm{def}}{=} (x_{CSC}^{s, \sigma}&(t; \tau_1, \tau_2), y_{CSC}^{s, \sigma}(t; \tau_1, \tau_2),\\
        &\varphi_{CSC}^{s, \sigma}(t; \tau_1, \tau_2)) \in \mathbb{R}^2 \times \mathbb{S}.
    \end{split}
\end{equation*}
Therefore, introducing domains of parameters $s$, $\sigma$, $\tau_1$, $\tau_2$ we obtain that
\begin{equation*}
    \begin{split}
        \mathcal{E}_{CSC}(t) = \{P_{CSC}^{s, \sigma}(t; \tau_1, \tau_2): \tau_1& \in [0, 2\pi), \tau_2 \in [\tau_1, t],\\
        &(s, \sigma) \in \mathbb{B}^2\}.
    \end{split}
\end{equation*}
This reasoning is valid also for the analytical description of $\mathcal{E}_{CCC}(t)$. Let $P^s_{CCC}(t; \tau_1, \tau_2) \in \mathcal{E}_{CCC}(t)$. Then the integrating of the system \eqref{eq:dynamics} using the initial conditions \eqref{eq:initial_conditions} and the control \eqref{eq:implicit_control_CCC} gives
\begin{equation}\label{eq:CCC_integrated}
    \begin{split}
        x^s_{CCC}(t; \tau_1, &\tau_2) = s(2\cos\tau_1 - 1 - 2\cos(2\tau_1 - \tau_2)\\
        &+ \cos(2\tau_1 - 2\tau_2 + t)),\\
        y^s_{CCC}(t; \tau_1, &\tau_2) = 2\sin\tau_1 - 2\sin(2\tau_1 - \tau_2)\\
        &+ \sin(2\tau_1 - 2\tau_2 + t),\\
        \varphi^s_{CCC}(t; \tau_1, &\tau_2) = \pi/2 + s(2\tau_1 - 2\tau_2 + t),
    \end{split}
\end{equation}
where we also suppose that
\begin{equation*}
    \begin{split}
        P^s_{CCC}(t; \tau_1, \tau_2) \overset{\mathrm{def}}{=} (x^s_{CCC}&(t; \tau_1, \tau_2), y^s_{CCC}(t; \tau_1, \tau_2),\\
        &\varphi^s_{CCC}(t; \tau_1, \tau_2)) \in \mathbb{R}^2 \times \mathbb{S}.
    \end{split}
\end{equation*}
Finally we obtain that
\begin{equation*}
    \begin{split}
        \mathcal{E}_{CCC}(t) = \{P^s_{CCC}(t; \tau_1, &\tau_2): \tau_1 \in [0, 2\pi), \tau_2 \in [\tau_1, t],\\
        &\tau_2 - \tau_1 \in [0, 2\pi), s \in \mathbb{B}\}.
    \end{split}
\end{equation*}

\subsection{Properties of $\mathcal{E}(t)$}
The analytical description of $\mathcal{E}(t)$ allows to state the following lemma.
\begin{lem}
    $\mathcal{E}(t)$ is a closed set at any moment of time $t \in \mathbb{R}^+_0$.
\end{lem}
\begin{pf}
    Since $P^{s, \sigma}_{CSC}(t; \tau_1, \tau_2)$, $P^s_{CCC}(t; \tau_1, \tau_2)$ depend continuously on $\tau_1$, $\tau_2$, only strict inequalities in the analytical descriptions of $\mathcal{E}_{CSC}(t)$, $\mathcal{E}_{CCC}(t)$ may break the closeness of $\mathcal{E}(t)$. Let's prove that $\mathcal{E}_{CSC}(t)$ is closed. A strict inequality is used only for $\tau_1 \in [0, 2\pi)$ in the explicit expression of $\mathcal{E}_{CSC}(t)$. Let $\tau_1 = 2\pi$, $\tau_2 \in [\tau_1, t]$, $t \geq 2\pi$. Then, calculations give 
    \begin{equation*}
        P^{s, \sigma}_{CSC}(t; 2\pi, \tau_2) = P^{s, \sigma}_{CSC}(t; 0, \tau_2 - 2\pi) \in \mathcal{E}_{CSC}(t).
    \end{equation*}
    Now, let's show that $\mathcal{E}_{CCC}(t)$ is closed. Strict inequalities are used only for $\tau_1 \in [0, 2\pi)$ and $\tau_2 - \tau_1 \in [0, 2\pi)$ in the analytical description of $\mathcal{E}_{CCC}(t)$. If $\tau_1 = 2\pi$, $\tau_2 \in [\tau_1, t]$, $\tau_2 - \tau_1 \in [0, 2\pi)$, $t \geq 2\pi$, then
    \begin{equation*}
        P^s_{CCC}(t; 2\pi, \tau_2) = P^s_{CCC}(t; 0, \tau_2 - 2\pi) \in \mathcal{E}_{CCC}(t).
    \end{equation*}
    If $\tau_1 = 2\pi$, $\tau_2 \in [\tau_1, t]$, $\tau_2 - \tau_1 = 2\pi$, $t \geq 4\pi$, then
    \begin{equation*}
        P^s_{CCC}(t; 2\pi, 4\pi) = P^s_{CCC}(t; 0, 0) \in \mathcal{E}_{CCC}(t).
    \end{equation*}
    If $\tau_1 \in [0, 2\pi)$, $\tau_2 \in [\tau_1, t]$, $\tau_2 - \tau_1 = 2\pi$, $t \geq \tau_1 + 2\pi$, then
    \begin{equation*}
        P^s_{CCC}(t; \tau_1, \tau_1 + 2\pi) = P^s_{CCC}(t; \tau_1, \tau_1) \in \mathcal{E}_{CCC}(t).
    \end{equation*}
\end{pf}

For the further purposes we formulate some key definitions related to closed multi-valued mappings.
\begin{defn}
    A closed multi-valued mapping $\mathcal{M}$ is said to be lower semicontinuous at the time moment $t_0$ when for any sufficiently small $\varepsilon > 0$ and any point $P_0 \in \mathcal{M}(t_0)$ there exists $\delta > 0$ such that for any $t: |t - t_0| < \delta$ there exists a point $P \in \mathcal{M}(t)$ such that $\rho(P, P_0) < \varepsilon$. 
\end{defn}
\begin{defn}
    A closed multi-valued mapping $\mathcal{M}$ is said to be upper semicontinuous at the time moment $t_0$ when for any sufficiently small $\varepsilon > 0$ there exists $\delta > 0$ such that for any $t: |t - t_0| < \delta$ and any point $P \in \mathcal{M}(t)$ there exists a point $P_0 \in \mathcal{M}(t_0)$ such that $\rho(P, P_0) < \varepsilon$.    
\end{defn}
\begin{defn}
    If a closed multi-valued mapping $\mathcal{M}$ is lower and upper semicontinuous at the time moment $t_0$, then it is continuous at time $t_0$.
\end{defn}

The continuity of the multi-valued mapping $\mathcal{R}$ is an established property \cite{lee1986foundations}, but a some multi-valued mapping $\mathcal{M}$ may not be continuous even if $\mathcal{M}(t) \subset \mathcal{R}(t)$ remains for any time $t \in \mathbb{R}^+_0$. For example, the evolution of the boundary of the Dubins car planar reachable set is not continuous at the some time moment \cite{buzikov2021time}.
\begin{lem}
    The multi-valued mapping $\mathcal{E}$ is continuous at any time moment $t_0 > 0$.
\end{lem}
\begin{pf}
    Let's prove that $\mathcal{E}$ is lower semicontinuous. For any point $P_0 = P^{s, \sigma}_{CSC}(t_0; \tau_1, \tau_2) \in \mathcal{E}_{CSC}(t_0)$ we suppose that if $\tau_2 < t_0$, then $P = P^{s, \sigma}_{CSC}(t; \tau_1, \tau_2)$, if $\tau_2 = t_0$ and $\tau_1 < \tau_2$, then $P = P^{s, \sigma}_{CSC}(t; \tau_1, t)$, if $\tau_2 = \tau_1 = t_0$, then $P = P^{s, \sigma}_{CSC}(t; t, t)$. Let's suppose $\delta = \varepsilon / \sqrt{2}$. Now, if $\tau_2 < t_0$ or $\tau_2 = \tau_1 = t_0$, then
    \begin{equation*}
        \begin{split}
            \rho(P, P_0) &\leq \sqrt{2 - 2\cos(t - t_0) + (t - t_0)^2}\\
            &\leq \sqrt{2}|t - t_0| < \delta\sqrt{2} = \varepsilon,
        \end{split}
    \end{equation*}
    and if $\tau_2 = t_0$ and $\tau_1 < \tau_2$, then
    \begin{equation*}
        \rho(P, P_0) \leq |t - t_0| < \delta < \varepsilon.
    \end{equation*}
    For any $P_0 = P^s_{CCC}(t_0; \tau_1, \tau_2) \in \mathcal{E}_{CCC}(t_0)$ we suppose that if $\tau_2 < t_0$, then $P = P^s_{CCC}(t; \tau_1, \tau_2)$; if $\tau_2 = t_0$ and $\tau_1 < \tau_2$, then $P = P^s_{CCC}(t; \tau_1, t)$; if $\tau_2 = \tau_1 = t_0$, then $P = P^s_{CCC}(t; t, t)$. Now, we have
    \begin{equation*}
        \begin{split}
            \rho(P, P_0) &\leq \sqrt{2 - 2\cos(t - t_0) + (t - t_0)^2}\\
            &\leq \sqrt{2}|t - t_0| < \delta\sqrt{2} = \varepsilon.
        \end{split}
    \end{equation*}
    Now, let's show that $\mathcal{E}$ is upper semicontinuous. Let's choose $\delta = \varepsilon / \sqrt2$. For any time $t: |t - t_0| < \delta$ and any point $P = P^{s, \sigma}_{CSC}(t; \tau_1, \tau_2) \in \mathcal{E}_{CSC}(t)$ we suppose that if $\tau_2 < t$, then $P_0 = P^{s, \sigma}_{CSC}(t_0; \tau_1, \tau_2)$; if $\tau_2 = t$ and $\tau_1 < \tau_2$, then $P_0 = P^{s, \sigma}_{CSC}(t_0; \tau_1, t_0)$; if $\tau_2 = \tau_1 = t$, then $P_0 = P^{s, \sigma}_{CSC}(t_0; t_0, t_0)$. For any point $P = P^s_{CCC}(t; \tau_1, \tau_2) \in \mathcal{E}_{CCC}(t)$ we suppose that if $\tau_2 < t$, then $P_0 = P^s_{CCC}(t_0; \tau_1, \tau_2)$; if $\tau_2 = t$ and $\tau_1 < \tau_2$, then $P_0 = P^s_{CCC}(t_0; \tau_1, t_0)$; if  $\tau_2 = \tau_1 = t$, then $P_0 = P^s_{CCC}(t_0; t_0, t_0)$. By the same way we conclude that $\rho(P, P_0) < \varepsilon$.
\end{pf}

\subsection{Location of the optimal interception point}
The listed above statements allow to prove the main theorem of this section. This theorem describes the configuration of the optimal interception point lying on the target trajectory.
\begin{thm}\label{th:optimal_interception_point}
    If the minimum-time lateral interception of the moving target by a Dubins car is possible, then the optimal interception time $T$ is finite, and the optimal interception point $E(T) \in \Tilde{\mathcal{B}}(T)$, where
    \begin{equation*}
        \Tilde{\mathcal{B}}(T) \overset{\mathrm{def}}{=} \{P' \in \mathcal{E}(T): \lim_{t \to T - 0} \min_{P \in \mathcal{B}(t)} \rho(P, P') = 0\}.
    \end{equation*}
\end{thm}
\begin{pf}
    Let $\mathcal{M}$ be an arbitrary continuous multi-valued mapping, such that $\mathcal{M}(t)$ is a closed set. Let
    \begin{equation*}
        \rho_\mathcal{M}(t) \overset{\mathrm{def}}{=} \min_{P \in \mathcal{M}(t)} \rho(P, E(t)).
    \end{equation*}
    Since the vector-function $E$ is continuous then $\rho_\mathcal{M}$ is a continuous function. Therefore, $\rho_\mathcal{R}$, $\rho_\mathcal{E}$ are continuous, but $\rho_\mathcal{B}$ is not necessarily continuous. For any time $t > 0$ the following inequalities hold
    \begin{equation*}
        \rho_\mathcal{R}(t) \leq \rho_\mathcal{E}(t) \leq \rho_\mathcal{B}(t),
    \end{equation*}
    since $\mathcal{B}(t) \subset \mathcal{E}(t) \subset \mathcal{R}(t)$. According to the definition of the minimum-time interception it holds at any time $t \in [0, T)$, that $E(t) \notin \mathcal{R}(t)$ and $E(T) \in \mathcal{R}(T)$. It implies that $\rho_{\mathcal{R}}(t) > 0$ at any time $t \in [0, T)$ and $\rho_{\mathcal{R}}(T) = 0$. Since $\mathcal{R}(t)$ is closed, $E(t) \notin \mathcal{R}(t)$, and $\mathcal{B}(t)$ is a boundary of $\mathcal{R}(t)$, we deduce that $\rho_\mathcal{R}(t) = \rho_\mathcal{B}(t)$ at any time $t \in [0, T)$. Now we conclude that
    \begin{equation*}
        \lim_{t \to T - 0}\rho_{\mathcal{E}}(t) \leq \lim_{t \to T - 0}\rho_{\mathcal{B}}(t) = \lim_{t \to T - 0}\rho_{\mathcal{R}}(t) = 0.
    \end{equation*}
    From the continuity of $\rho_\mathcal{E}$ we get $\rho_{\mathcal{E}}(T) = 0$, so $E(T) \in \mathcal{E}(T)$. According to the $\Tilde{\mathcal{B}}(T)$ definition we obtain that $E(T) \in \Tilde{\mathcal{B}}(T)$, since $\rho_{\mathcal{B}}(T-0) = 0$.
\end{pf}

Thus the analytical description of $\mathcal{E}(T) \supset \Tilde{\mathcal{B}}(T)$ gives the desired algebraic equations which solution allows to determine the optimal control parameters. But it is still possible to reduce the space of optimal control candidates. Consider the next mapping for the following purposes:
\begin{equation*}
    S_\tau(x, y, \varphi) = (x + \tau\cos\varphi, y + \tau\sin\varphi, \varphi).
\end{equation*}
Note that if the Dubins car is located at the point $P_0 = (x_0, y_0, \varphi_0) \in \mathcal{R}(t_0)$ at the time moment $t_0$ and its control $u(t) = 0$ for any $t \geq t_0$, then the Dubins car is located at the point $S_{t - t_0}(P_0) \in \mathcal{R}(t)$
at time $t$.

\begin{lem}\label{lem:internal_point}
    If a point $P_0 = (x_0, y_0, \varphi_0)$ is an internal point of $\mathcal{R}(t_0)$, then for any time moment $t > t_0$ there exists $\delta_0 > 0$ such that for any $\delta \in [0, \delta_0)$ the point $S_{t - t_0}(P_0)$ is an internal point of $\mathcal{R}(t-\delta)$.
\end{lem}

\begin{pf}
    Since $P_0$ is an internal point, there exists a ball $B_{\delta_0}(P_0) = \{P\in\mathbb{R}^2 \times \mathbb{S}: \rho(P, P_0) < \delta_0\}$ which is a subset of $\mathcal{R}(t_0)$. Let's demand that $\delta_0 < t - t_0$ and consider an arbitrary $\delta \in [0, \delta_0)$. According to properties of the mapping $S_{t-t_0-\delta}$:
    \begin{equation*}
        S_{t-t_0-\delta}(B_{\delta_0}(P_0)) \subset \mathcal{R}(t-\delta).
    \end{equation*}
    It is obvious that $S_\delta(P_0) \in B_{\delta_0}(P_0)$ since
    \begin{equation*}
        \rho(S_\delta(P_0), P_0) = \delta < \delta_0.
    \end{equation*}
    Moreover, $S_\delta(P_0)$ is an internal point of $B_{\delta_0}(P_0)$. The mapping $S_{t - t_0 - \delta}$ is a continuous mapping. Hence, an image $S_{t - t_0 - \delta}(B_{\delta_0}(P_0)) \subset \mathcal{R}(t - \delta)$ is an open set and the point $S_{t-t_0-\delta}(S_\delta(P_0)) = S_{t-t_0}(P_0)$ is an internal point of the image. 
\end{pf}

\begin{lem}\label{lem:exclude_some_cycled_trajectories}
    Let $t - \tau_2 \geq 2\pi$ and $\tau_2 - \tau_1 > 0$. If one of the inequalities  $\tau_1 > 0$ or $t - \tau_2 > 2\pi$ is fulfilled, then $P^{s, \sigma}_{CSC}(t; \tau_1, \tau_2) \notin \Tilde{\mathcal{B}}(t)$.
\end{lem}

\begin{pf}
    According to $\Tilde{\mathcal{B}}(t)$ definition it is sufficient to prove that there exists $\delta_0$ such that for any $\delta \in [0, \delta_0)$ the point $P^{s, \sigma}_{CSC}(t; \tau_1, \tau_2)$ holds an internal point for $\mathcal{R}(t-\delta)$. Suppose that $\tau = (\tau_2 - \tau_1) / 2$ and consider a point
    \begin{equation*}
        P^{s, \sigma}_{CSC}(t-\tau; \tau_1, \tau_2 - \tau) \in \mathcal{R}(t - \tau).
    \end{equation*}
    This point can be attained using the trajectory that contains a cycle. This cycle could be shifted to any part of the trajectory and according to Lemma~4 from~\cite{patsko2003three} the point $P^{s, \sigma}_{CSC}(t-\tau; \tau_1, \tau_2 - \tau)$ is an internal point of $\mathcal{R}(t - \tau)$. Hence, according to Lemma~\ref{lem:internal_point} we conclude that $P^{s, \sigma}_{CSC}(t; \tau_1, \tau_2)$ holds an internal point for $\mathcal{R}(t-\delta)$ since
    \begin{equation*}
        S_\tau(P^{s, \sigma}_{CSC}(t-\tau; \tau_1, \tau_2 - \tau)) = P^{s, \sigma}_{CSC}(t; \tau_1, \tau_2).
    \end{equation*}
\end{pf}

The following lemma allows to avoid a case with two cycles when we consider an interception problem.
\begin{lem}\label{lem:4pi_inner_point}
    The configuration $(0, 0, \pi/2) \notin \Tilde{\mathcal{B}}(4\pi)$.
\end{lem}

\begin{pf}
    According to the definition of $\Tilde{\mathcal{B}}$ it is enough to prove that there exists $\delta_0 > 0$ such that for any $\delta \in [0, \delta_0)$ the point $(0, 0, \pi/2)$ is an internal point of $\mathcal{R}(4\pi - \delta)$. Let's consider the control given by 
    \begin{equation*}
        u(t) = 
        \begin{cases}
            0,& \quad t \in [0, \pi/2 + \chi),\\
            2/3 + \xi,& \quad t \in [\pi/2 + \chi, 2\pi + \chi),\\
            0,& \quad t \in [2\pi + \chi, 5\pi/2 + \chi),\\
            2/3 + \eta,& \quad t \in [5\pi/2 + \chi, +\infty).\\
        \end{cases}
    \end{equation*}
    Integration of the system \eqref{eq:dynamics} with this control function and the initial conditions given by \eqref{eq:initial_conditions} leads to
    \begin{equation*}
        \begin{split}
            x = &\frac{(\eta - \xi)\cos(\pi(1 + \frac{3\xi}2))}{(\frac23 + \xi)(\frac23 + \eta)} - \frac\pi2\sin\left(\pi\left(1 + \frac{3\xi}2\right)\right)\\
            & + \frac{\cos(\pi(1 + \frac{3\xi}2) + (\frac23 + \eta)(t - \frac{5\pi}2 - \chi))}{\frac23+\eta} - \frac1{\frac23 + \xi},\\
        \end{split}
    \end{equation*}
    \begin{equation*}
        \begin{split}
            y = &\frac{(\eta - \xi)\sin(\pi(1 + \frac{3\xi}2))}{(\frac23 + \xi)(\frac23 + \eta)} + \frac\pi2\cos\left(\pi\left(1 + \frac{3\xi}2\right)\right)\\
            & + \frac{\sin(\pi(1 + \frac{3\xi}2) + \left(\frac23 + \eta\right)(t - \frac{5\pi}2 - \chi))}{\frac23+\eta} + \frac\pi2 + \chi,\\
        \end{split}
    \end{equation*}
    \begin{equation*}
        \varphi = \frac{3\pi}2(1 + \xi) + \left(\frac23 + \eta\right)\left(t - \frac{5\pi}2 - \chi\right)
    \end{equation*}
    at time $t \geq 5\pi/2 + \chi$. If $t = 4\pi$, $\xi = 0$, $\eta = 0$, $\chi = 0$, then the point $(x, y, \varphi) = (0, 0, 5\pi/2) = (0, 0, \pi/2) \in \mathcal{R}(4\pi)$. The calculation of the Jacobi matrix determinant leads to the following value 
    \begin{equation*}
        \left|\frac{\partial(x, y, \varphi)}{\partial(\xi, \eta, \chi)}\right|_{\substack{\xi = 0 \\ \eta = 0 \\ \chi = 0 \\ t = 4\pi}} =
        \begin{vmatrix}
            \frac{18+3\pi^2}4 & -\frac92 & 0\\
            \frac{9\pi}2 & \frac{9\pi}2 & 0\\
            \frac{3\pi}2 & \frac{3\pi}2 & -\frac23
        \end{vmatrix}
        = -\frac{27\pi}2 - \frac{9\pi^3}8 \neq 0.
    \end{equation*}
    Since this value does not vanish and the corresponding mapping of $(\xi, \eta, \chi)$ to $(x, y, \varphi)$ is continuous for small enough $\delta_0 > 0$ at any time $t \in (4\pi - \delta_0, 4\pi + \delta_0)$, this mapping transforms any small enough neighborhood of $(\xi, \eta, \chi) = (0, 0, 0)$ to the neighborhood of $(x, y, \varphi) = (0, 0, \pi/2)$. The constructed control is admissible and determined completely by $(\xi, \eta, \chi)$. It implies that the neighborhood of $(x, y, \varphi) = (0, 0, \pi/2)$ is a subset of $\mathcal{R}(t)$ at any time $t \in  (4\pi - \delta_0, 4\pi + \delta_0)$. Therefore, the point $(0, 0, \pi/2)$ is an internal point of $\mathcal{R}(4\pi-\delta)$ for any $\delta \in [0, \delta_0)$.
\end{pf}

Now it is possible to prove a statement which is a generalization of \cite[Lemma 2]{techy2009minimum}. According to Lemma \ref{lem:4pi_inner_point} and Theorem~\ref{th:optimal_interception_point} the point $(0, 0, \pi/2)$ is not an optimal interception point at time $t = 4\pi$. Moreover, Lemma~\ref{lem:exclude_some_cycled_trajectories} excludes some trajectories with cycles. We use it for the optimal control candidates reduction in the following theorem.

\begin{thm}\label{th:CSC_CCC_additional_restriction}
    For attaining any point $P_0 \in \Tilde{\mathcal{B}}(t_0)$ at time $t_0 \in \mathbb{R}^+_0$ it suffices to use the controls \eqref{eq:implicit_control_CSC}, \eqref{eq:implicit_control_CCC} with $t_0 - \tau_2 \in [0, 2\pi)$ or the controls \eqref{eq:implicit_control_CSC}, \eqref{eq:implicit_control_CCC} with $\tau_1 = 0$, $t_0 - \tau_2 = 2\pi$.
\end{thm}

\begin{pf}
    Suppose that the controls \eqref{eq:implicit_control_CSC}, \eqref{eq:implicit_control_CCC} with $t_0 - \tau_2 \geq 4\pi$ lead to the point $P_0 \in \Tilde{\mathcal{B}}(t_0)$. According to the Bellman's principle these controls remain optimal from the time $t_0 - 4\pi$ to $t_0$. Considering the new time-optimal interception problem from the time $t_0 - 4\pi$ raises a contradiction with Lemma~\ref{lem:4pi_inner_point} and Theorem~\ref{th:optimal_interception_point} since the new problem interception point is located at $(0, 0, \pi/2)$ at time $t = 4\pi$. So, the controls \eqref{eq:implicit_control_CSC}, \eqref{eq:implicit_control_CCC} with $t_0 - \tau_2 \geq 4\pi$ are redundant.
    
    Now, let's suppose that the control \eqref{eq:implicit_control_CSC} with $t_0 - \tau_2 \geq 2\pi$ leads to the point $P_0 = P^{s, \sigma}_{CSC}(t_0; \tau_1, \tau_2) \in \Tilde{\mathcal{B}}(t_0)$. We associate the case when $\tau_2 - \tau_1 = 0$ with CCC-trajectories and consider it further. If Lemma~\ref{lem:exclude_some_cycled_trajectories} cases hold, then $P_0 \notin \Tilde{\mathcal{B}}(t_0)$. Otherwise, $\tau_1 = 0$ and $t_0 - \tau_2 = 2\pi$.
    
    Next, let's suppose that the control \eqref{eq:implicit_control_CCC} with $t_0 - \tau_2 \geq 2\pi$ leads to the point $P_0 = P^s_{CCC}(t_0; \tau_1, \tau_2) \in \Tilde{\mathcal{B}}(t_0)$. Consider the case when $\tau_1 > 0$. According to \cite[Lemma 2]{patsko2003three} if
    \begin{equation*}
        \tau_1 + (t_0 - \tau_2) > \tau_2 - \tau_1,
    \end{equation*}
    then $P_0$ is an internal point of $\mathcal{R}(t_0)$. Consider this lemma for a terminal point at time $t \in (t_0 - \tau_1 / 2, t_0]$. It follows that the terminal point holds internal for $\mathcal{R}(t)$ since
    \begin{equation*}
        \tau_1 + t - \tau_2 > \tau_1 + t_0 - \frac{\tau_1}2 - \tau_2 = \frac{\tau_1}2 + t_0 - \tau_2 > 2\pi > \tau_2 - \tau_1.
    \end{equation*}
    Therefore, $P_0 \notin \Tilde{\mathcal{B}}(t_0)$. Now suppose that $\tau_1 = 0$ and $\tau_2 - \tau_1 > 0$. Let's consider a new control $u^{-s}_{CCC}(\cdot; \tau_1', \tau_2')$ with $\tau_1' = \tau_2$, $\tau_2' - \tau_1' = t_0 - \tau_2 - 2\pi$ which leads to the same point $P_0$. For the arbitrary $t \in (t_0 - \tau_1'/2, t_0]$ the corresponding terminal point also holds internal for $\mathcal{R}(t)$ since
    \begin{equation*}
        \tau_1' + t - \tau_2' > \frac{\tau_1'}2 + t_0 - \tau_2' = \frac{\tau_1'}2 + 2\pi > 2\pi > \tau_2' - \tau_1'.
    \end{equation*}
    If $\tau_1 = 0$ and $\tau_2 - \tau_1 = 0$, then the control $u^{-s}_{CCC}(\cdot; \tau_1', \tau_2')$ with $\tau_1' = 0$, $t_0 - \tau_2' = 2\pi$ leads to the point $P_0$.
\end{pf}

\section{Optimal control synthesis}\label{sec:optimal_control}
In this section, we consider some scenarios of the optimal interception. These scenarios are enough to completely solve the problem. Finally, an algorithm of the optimal control selection is described.

\subsection{Optimal CSC-trajectory}
According to Theorem~\ref{th:optimal_interception_point} the optimal interception point $E(T) \in \Tilde{\mathcal{B}}(T) \subset \mathcal{E}(T)$. Suppose that the using of the CSC-trajectory is optimal. It follows from \eqref{eq:terminal_conditions} that
\begin{equation*}
    E(T) = P^{s, \sigma}_{CSC}(T; \tau_1, \tau_2).
\end{equation*}
Next, we use the equations given by \eqref{eq:CSC_integrated} and infer that
\begin{equation}\label{eq:CSC_interception}
    \begin{split}
        \xi^{s, \sigma}_{CSC}(T) \overset{\mathrm{def}}{=}& sx_E(T) + 1 - s\sigma\sin\varphi_E(T)\\
        &= (1 - s\sigma)\cos\tau_1 - (\tau_2 - \tau_1)\sin\tau_1,\\
        \eta^\sigma_{CSC}(T) \overset{\mathrm{def}}{=}& y_E(T) + \sigma\cos\varphi_E(T)\\
        &= (1 - s\sigma)\sin\tau_1 + (\tau_2 - \tau_1)\cos\tau_1,\\
        \varphi_E(T) = \frac\pi2& + s\tau_1 + \sigma(T - \tau_2) + 2\pi k, \quad k \in \mathbb{Z}.
    \end{split}
\end{equation}
Consider a new function given by
\begin{equation*}
    \begin{split}
         \rho^{2, s, \sigma}_{CSC}(T) \overset{\mathrm{def}}{=} \left(\xi^{s, \sigma}_{CSC}(T)\right)^2 + &\left(\eta^\sigma_{CSC}(T)\right)^2\\
         &= (1 - s\sigma)^2 + (\tau_2 - \tau_1)^2.
    \end{split}
\end{equation*}
According to Lemma~\ref{lem:reduced_CSC_CCC} we suppose that $\tau_2 \geq \tau_1$. It gives the following equation
\begin{equation*}
    \tau_2 - \tau_1 = \sqrt{\rho^{2, s, \sigma}_{CSC}(T) - (1 - s\sigma)^2}.
\end{equation*}
Now we substitute the value $\tau_2 - \tau_1$ to the first and second equations of \eqref{eq:CSC_interception} and after some algebra it provides the following
\begin{equation*}
    \begin{split}
        C^{s, \sigma}_{CSC}(T) \overset{\mathrm{def}}{=} \frac{\eta^\sigma_{CSC}(T)}{\rho^{2, s, \sigma}_{CSC}(T)}&\sqrt{\rho^{2, s, \sigma}_{CSC}(T) - (1 - s\sigma)^2}\\
        &+ \frac{(1 - s\sigma)\xi^{s, \sigma}_{CSC}(T)}{\rho^{2, s, \sigma}_{CSC}(T)} = \cos\tau_1,\\
        S^{s, \sigma}_{CSC}(T) \overset{\mathrm{def}}{=} -\frac{\xi^{s, \sigma}_{CSC}(T)}{\rho^{2, s, \sigma}_{CSC}(T)}&\sqrt{\rho^{2, s, \sigma}_{CSC}(T) - (1 - s\sigma)^2}\\
        &+ \frac{(1 - s\sigma)\eta^\sigma_{CSC}(T)}{\rho^{2, s, \sigma}_{CSC}(T)} = \sin\tau_1.
    \end{split}
\end{equation*}
For the further purposes we use the next function:
\begin{equation*}
    \arctan_2(y, x) \overset{\mathrm{def}}{=}
        \begin{cases}
            &\arccos{\cfrac{x}{\sqrt{x^2 + y^2}}}, \quad y \geq 0,\\
            &2\pi - \arccos{\cfrac{x}{\sqrt{x^2 + y^2}}}, \quad y < 0.
        \end{cases}
\end{equation*}
The previous equations and $\tau_1 \in [0, 2\pi)$ provide the values
\begin{equation*}
    \begin{split}
        &\tau_1 = \theta^{1, s, \sigma}_{CSC}(T) \overset{\mathrm{def}}{=} \arctan_2\left(S^{s, \sigma}_{CSC}(T), C^{s, \sigma}_{CSC}(T)\right),\\
        &\tau_2 = \theta^{2, s, \sigma}_{CSC}(T) \overset{\mathrm{def}}{=} \theta^{1, s, \sigma}_{CSC}(T) + \sqrt{\rho^{2, s, \sigma}_{CSC}(T) - (1 - s\sigma)^2}.
    \end{split}
\end{equation*}
The third equation of \eqref{eq:CSC_interception} and $T \geq \tau_2$ give the following problem:
\begin{equation*}
    \left\{\begin{aligned}
        &\varphi_E(T) - \frac\pi2 - s\theta^{1, s, \sigma}_{CSC}(T) - \sigma(T - \theta^{2, s, \sigma}_{CSC}(T)) = 2\pi k,\\
        &T \geq \theta^{2, s, \sigma}_{CSC}(T), \quad k \in \mathbb{Z}.
    \end{aligned}\right.
\end{equation*}
The additional constraint from the Theorem~\ref{th:CSC_CCC_additional_restriction} restricts the value of $T - \tau_2 < 2\pi$ in the cycle-free case. Therefore, the optimal interception time $T$ is subject to
\begin{equation}\label{eq:optimal_T_for_CSC}
    \begin{split}
        F^{s, \sigma}_{CSC}&(T) \overset{\mathrm{def}}{=} -T + \theta^{2, s, \sigma}_{CSC}(T)\\
        &+ \mathrm{mod}\left(\sigma\left(\varphi_E(T) - \frac\pi2\right)-s\sigma\theta^{1, s, \sigma}_{CSC}(T), 2\pi\right) = 0.
    \end{split}
\end{equation}
Hence, we have four equations \eqref{eq:optimal_T_for_CSC} since the trajectory configuration $(s, \sigma) \in \mathbb{B}^2$ is not established. For each configuration $(s, \sigma)$ we should find a minimal non-negative solution $T^{s, \sigma}_{CSC}$ of the algebraic equation \eqref{eq:optimal_T_for_CSC} which is transcendental in general case. If this equation does not have a solution for some configuration $(s, \sigma) \in \mathbb{B}^2$, then we suppose that $T^{s, \sigma}_{CSC} = +\infty$.
The optimal control $u^{s, \sigma}_{CSC}(\cdot; \tau_1, \tau_2)$ is defined by \eqref{eq:implicit_control_CSC} with the following choice of parameters:
\begin{equation}\label{eq:CSC_optimal_params}
    \begin{split}
        &(s, \sigma) = \arg\min_{(\tilde{s}, \tilde{\sigma}) \in \mathbb{B}^2} T^{\tilde{s}, \tilde{\sigma}}_{CSC},\\
        &\tau_1 = \theta^{1, s, \sigma}_{CSC}\left(T^{s, \sigma}_{CSC}\right), \quad \tau_2 = \theta^{2, s, \sigma}_{CSC}\left(T^{s, \sigma}_{CSC}\right).
    \end{split}
\end{equation}
For the optimal configuration $(s, \sigma)$ we suppose that the optimal interception time is $T_{CSC} = T^{s, \sigma}_{CSC}$ in the considered case of the cycle-free CSC-trajectories.

If there is a cycle in the optimal CSC-trajectory and a non-vanished S-part, then according to Theorem~\ref{th:CSC_CCC_additional_restriction} $\tau_1 = 0$, $T - \tau_2 = 2\pi$. In this case the optimal interception time $T$ is subject to
\begin{equation}\label{eq:optimal_T_for_CSC_with_cycle}
    F_{SC}(T) \overset{\mathrm{def}}{=} \rho(E(T), (0, T-2\pi, \pi/2)) = 0.
\end{equation}
We should find a minimal solution $T_{SC} \geq 2\pi$ of the algebraic equation \eqref{eq:optimal_T_for_CSC_with_cycle} which is transcendental in general case. If this equation does not have a solution $T_{SC} \geq 2\pi$, then we suppose that the optimal interception time $T_{SC} = +\infty$. The optimal control $u^{s, \sigma}_{CSC}(\cdot; \tau_1, \tau_2)$ is defined by \eqref{eq:implicit_control_CSC} with the arbitrary choice of parameters $(s, \sigma) \in \mathbb{B}^2$ and $\tau_1 = 0$, $\tau_2 = T_{SC} - 2\pi$.

If the S-part in the optimal CSC-trajectory vanishes, then we relate this case to the CCC-trajectory case which is considered further.

\subsection{Optimal CCC-trajectory}
Let's consider that the using of the CCC-trajectory is optimal. The line of reasoning is similar to the CSC-trajectory case. It follows from \eqref{eq:terminal_conditions} that
\begin{equation*}
    E(T) = P^s_{CCC}(T; \tau_1, \tau_2).
\end{equation*}
Next, we use the equations given by \eqref{eq:CCC_integrated} and infer that
\begin{equation}\label{eq:CCC_interception}
    \begin{split}
        \xi^s_{CCC}(T&) \overset{\mathrm{def}}{=} \left(sx_E(T) + 1 - \sin\varphi_E(T)\right)/2\\
        &= (1 - \cos(\tau_2 - \tau_1))\cos\tau_1 - \sin(\tau_2 - \tau_1)\sin\tau_1,\\
        \eta^s_{CCC}(T&) \overset{\mathrm{def}}{=} \left(y_E(T) + s\cos\varphi_E(T)\right)/2\\
        &= (1 - \cos(\tau_2 - \tau_1))\sin\tau_1 + \sin(\tau_2 - \tau_1)\cos\tau_1,\\ 
        \varphi_E(T) &= \pi/2 + s(2\tau_1 - 2\tau_2 + T) + 2\pi k, \quad k \in \mathbb{Z}.
    \end{split}
\end{equation}
Consider a new function given by
\begin{equation*}
    \begin{split}
         \rho^{2, s}_{CCC}(T) \overset{\mathrm{def}}{=} \left(\xi^s_{CCC}(T)\right)^2& + \left(\eta^s_{CCC}(T)\right)^2\\
         &= 2 - 2 \cos(\tau_2 - \tau_1).
    \end{split}
\end{equation*}
It follows that
\begin{equation*}
    \begin{split}
        &\cos(\tau_2 - \tau_1) = 1 - \frac{\rho^{2, s}_{CCC}(T)}2,\\
        &\sin(\tau_2 - \tau_1) =\mu\sqrt{1 - \left(1 - \frac{\rho^{2, s}_{CCC}(T)}2\right)^2},
    \end{split}
\end{equation*}
where $\mu \in \mathbb{B}$. According to Lemma~\ref{lem:reduced_CSC_CCC} we suppose that $\tau_2 - \tau_1 \in [0, 2\pi)$. If $\mu = 1$ then the length of the second circle arc of the CCC-trajectory is $\tau_2 - \tau_1 \in [0, \pi]$. Otherwise, if $\mu = -1$ then $\tau_2 - \tau_1 \in (\pi, 2\pi)$. Therefore,
\begin{equation*}
    \tau_2 = \tau_1 - \pi(\mu - 1) + \mu\arccos\left(1 - \frac{\rho^{2, s}_{CCC}(T)}2\right).
\end{equation*}
Now we substitute the values $\cos(\tau_2 - \tau_1)$ and $\sin(\tau_2 - \tau_1)$ to the first and second equations of \eqref{eq:CCC_interception} and after some algebra it provides the following
\begin{equation*}
    \begin{split}
        C^{s, \mu}_{CCC}(T) \overset{\mathrm{def}}{=} \mu\frac{\eta^s_{CCC}(T)}{\rho^{2, s}_{CCC}(T)}&\sqrt{1 - \left(1 - \frac{\rho^{2, s}_{CCC}(T)}2\right)^2}\\
        &+ \frac{\xi^s_{CCC}(T)}2 = \cos\tau_1,\\
        S^{s, \mu}_{CCC}(T) \overset{\mathrm{def}}{=} -\mu\frac{\xi^s_{CCC}(T)}{\rho^{2, s}_{CCC}(T)}&\sqrt{1 - \left(1 - \frac{\rho^{2, s}_{CCC}(T)}2\right)^2}\\
        &+ \frac{\eta^s_{CCC}(T)}2 = \sin\tau_1.\\
    \end{split}
\end{equation*}
The previous equations and $\tau_1 \in [0, 2\pi)$ provide the values
\begin{equation*}
    \begin{split}
        \tau_1 = \theta^{1, s, \mu}_{CCC}(T) \overset{\mathrm{def}}{=} &\arctan_2\left(S^{s, \mu}_{CCC}(T), C^{s, \mu}_{CCC}(T)\right),\\
        \tau_2 = \theta^{2, s, \mu}_{CCC}(T) \overset{\mathrm{def}}{=} &\theta^{1, s, \mu}_{CCC}(T) - \pi(\mu - 1)\\
        &+\:\mu\arccos\left(1 - \frac{\rho^{2, s}_{CCC}(T)}2\right).
    \end{split}
\end{equation*}
The third equation of \eqref{eq:CCC_interception} and $T \geq \tau_2$ give the following problem:
\begin{equation*}
    \left\{\begin{aligned}
        &\varphi_E(T) - \frac\pi2 - 2s\theta^{1, s, \mu}_{CCC}(T) + 2s\theta^{2, s, \mu}_{CCC}(T) - sT = 2\pi k,\\
        &T \geq \theta^{2, s, \mu}_{CCC}(T), \quad k\in \mathbb{Z}.
    \end{aligned}\right.
\end{equation*}
The additional constraint from the Theorem~\ref{th:CSC_CCC_additional_restriction} restricts the value of $T - \tau_2 < 2\pi$ in the cycle-free case. Therefore, the optimal interception time $T$ is subject to
\begin{equation}\label{eq:optimal_T_for_CCC}
    \begin{split}
        &F^{s, \mu}_{CCC}(T) \overset{\mathrm{def}}{=} -T + \theta^{2, s, \mu}_{CCC}(T)\:+ \\
        &\mathrm{mod}(s(\varphi_E(T) - \frac\pi2) - 2\theta^{1, s, \mu}_{CCC}(T) + \theta^{2, s, \mu}_{CCC}(T), 2\pi) = 0.
    \end{split}
\end{equation}
Hence, we have four equations \eqref{eq:optimal_T_for_CCC} since the trajectory configuration $(s, \mu) \in \mathbb{B}^2$ is not established. For each configuration $(s, \mu)$ we should find a minimal non-negative solution $T^{s, \mu}_{CCC}$ of the algebraic equation \eqref{eq:optimal_T_for_CCC} which is transcendental in general case. If this equation does not have a solution for some configuration $(s, \mu) \in \mathbb{B}^2$, then we suppose that $T^{s, \mu}_{CSC} = +\infty$. The optimal control $u^s_{CCC}(\cdot; \tau_1, \tau_2)$ is defined by \eqref{eq:implicit_control_CCC} with the following choice of parameters:
\begin{equation}\label{eq:CCC_optimal_params}
    \begin{split}
        &(s, \mu) = \arg\min_{(\tilde{s}, \tilde{\mu}) \in \mathbb{B}^2} T^{\tilde{s}, \tilde{\mu}}_{CCC},\\
        &\tau_1 = \theta^{1, s, \mu}_{CCC}\left(T^{s, \mu}_{CCC}\right), \quad \tau_2 = \theta^{2, s, \mu}_{CCC}\left(T^{s, \mu}_{CCC}\right).
    \end{split}
\end{equation}
For the optimal configuration $(s, \mu)$ we suppose that the optimal interception time is $T_{CCC} = T^{s, \mu}_{CCC}$ in the considered case of the cycle-free CCC-trajectories.

If there is a cycle in the optimal CCC-trajectory, then according to Theorem~\ref{th:CSC_CCC_additional_restriction} the first circle arc is vanished $\tau_1 = 0$ and $T - \tau_2 = 2\pi$. In this case, the optimal interception time $T$ is subject to
\begin{equation}\label{eq:optimal_T_for_CCC_with_cycle}
    \begin{split}
        F_{CC}(T) \overset{\mathrm{def}}{=} \rho(E(T), (&(1-\cos T)\:\mathrm{sgn}\:x_E(T), \sin T,\\
        &\pi/2 - T\:\mathrm{sgn}\:x_E(T))) = 0.
    \end{split}
\end{equation}
We should find a minimal solution $T_{CC} \geq 2\pi$ of the algebraic equation \eqref{eq:optimal_T_for_CCC_with_cycle} which is transcendental in general case. If this equation does not have a solution $T_{CC} \geq 2\pi$, then we suppose that the optimal interception time $T_{CC} = +\infty$. The optimal control $u^s_{CCC}(\cdot; \tau_1, \tau_2)$ is defined by \eqref{eq:implicit_control_CCC} with the parameters $s = \:\mathrm{sgn}\:x_E(T)$, $\tau_1 = 0$, $\tau_2 = T_{CC} - 2\pi$.

\subsection{Optimal control}

The definition of $\Tilde{\mathcal{B}}(T)$ allows to consider that any point of $\Tilde{\mathcal{B}}(T)$ is attainable using either CSC- or CCC-trajectories. We consider the brute force algorithm of searching of the optimal control. This algorithm calculates minimal roots of 10 equations and compares them.

It supposes that the non-negative minimal root $T^{s, \sigma}_{CSC}$ of the equation \eqref{eq:optimal_T_for_CSC} is defined by some root-finding algorithm for each $(s, \sigma) \in \mathbb{B}^2$. The corresponding optimal cycle-free CSC-trajectory is defined by the parameters \eqref{eq:CSC_optimal_params} and the interception time is $T_{CSC}$. If $T_{CSC} = +\infty$, then the equation \eqref{eq:optimal_T_for_CSC} does not have non-negative roots for each $(s, \sigma) \in \mathbb{B}^2$. It means that the minimum time interception along any cycle-free CSC-trajectory is impossible. The CSC-trajectory with a cycle and non-zero S-part is examined by solution $T_{SC} \geq 2\pi$ of the algebraic equation \eqref{eq:optimal_T_for_CSC_with_cycle}.

The algorithm also supposes that the non-negative minimal root $T^{s, \mu}_{CCC}$ of the equation \eqref{eq:optimal_T_for_CCC} is defined by some root-finding algorithm for each $(s, \mu) \in \mathbb{B}^2$. The corresponding optimal cycle-free CCC-trajectory is defined by the parameters \eqref{eq:CCC_optimal_params} and the interception time is $T_{CCC}$. If $T_{CCC} = +\infty$, then the equation \eqref{eq:optimal_T_for_CCC} does not have non-negative roots for each $(s, \mu) \in \mathbb{B}^2$. The cycled CCC-trajectory is examined by solution $T_{CC} \geq 2\pi$ of the algebraic equation \eqref{eq:optimal_T_for_CCC_with_cycle}. The final step of the algorithm is comparison of $T_{CSC}$, $T_{SC}$, $T_{CCC}$ and $T_{CC}$ values. If at least one of these values is finite then lateral interception is possible. The following theorem determines the brute force algorithm.

\begin{thm}\label{th:optimal_control}
    Let the minimum-time lateral interception of the moving target by the Dubins car be possible and $T^* = \min(T_{CSC}, T_{SC}, T_{CCC}, T_{CC})$. If $T^* = T_{CSC}$, then $u = u^{s, \sigma}_{CSC}(\cdot; \tau_1, \tau_2)$ is an optimal control, where $s$, $\sigma$, $\tau_1$, $\tau_2$ are defined by \eqref{eq:CSC_optimal_params}. If $T^* = T_{SC}$, then $u = u^{s, \sigma}_{CSC}(\cdot; \tau_1, \tau_2)$ is an optimal control, where $(s, \sigma) \in \mathbb{B}^2$, $\tau_1 = 0$, $\tau_2 = T_{SC} - 2\pi$. If $T^* = T_{CCC}$, then $u = u^s_{CCC}(\cdot; \tau_1, \tau_2)$ is an optimal control, where $s$, $\tau_1$, $\tau_2$ are defined by \eqref{eq:CCC_optimal_params}. If $T^* = T_{CC}$, then $u = u^s_{CCC}(\cdot; \tau_1, \tau_2)$ is an optimal control, where $s =\:\mathrm{sgn}\:x_E(T_{CC})$, $\tau_1 = 0$, $\tau_2 = T_{CC} - 2\pi$.
\end{thm}

\section{Examples}\label{sec:examples}
This section provides a set of optimal interception examples. The main goal of these examples is to highlight the necessity of considering all 10 algebraic equations in general case. Also, we demonstrate an application of the proposed brute force algorithm.

The Dubins car reachable set has a plane of symmetry $x = 0$ in $\mathbb{R}^2 \times \mathbb{S}$. For the CSC-optimal trajectory case it means that if there exists an example when the optimal interception of the target $E(t) = (x_E(t), y_E(t), \varphi_E(t))$ proceeds using the CSC-trajectory and the control $u^{s, \sigma}_{CSC}(\cdot; \tau_1, \tau_2)$ is optimal, then the other target
\begin{equation}\label{eq:mirror_target}
    E^*(t) = (-x_E(t), y_E(t), \pi - \varphi_E(t))
\end{equation}
could be optimally intercepted using the mirror control $u^{-s, -\sigma}_{CSC}(\cdot; \tau_1, \tau_2)$. The same situation occurs when the control $u^s_{CCC}(\cdot; \tau_1, \tau_2)$ is optimal. The target $E^*$ could be optimally intercepted using the control $u^{-s}_{CCC}(\cdot; \tau_1, \tau_2)$.

This symmetry allows considering only 6 cases for the demonstration that it is necessary to solve each of 10 equations in general case.

Here and further in figures, the target orientation on its Cartesian plane path is presented by the light triangles and Dubins car velocity direction by filled black ones. At the interception moment, the black triangular coincides with the light one. The continuous line demonstrates an optimal interception path on the plane.

\subsection{CSC-optimal examples}

\begin{figure}[h!]
    \begin{center}
        \includegraphics[width=\linewidth]{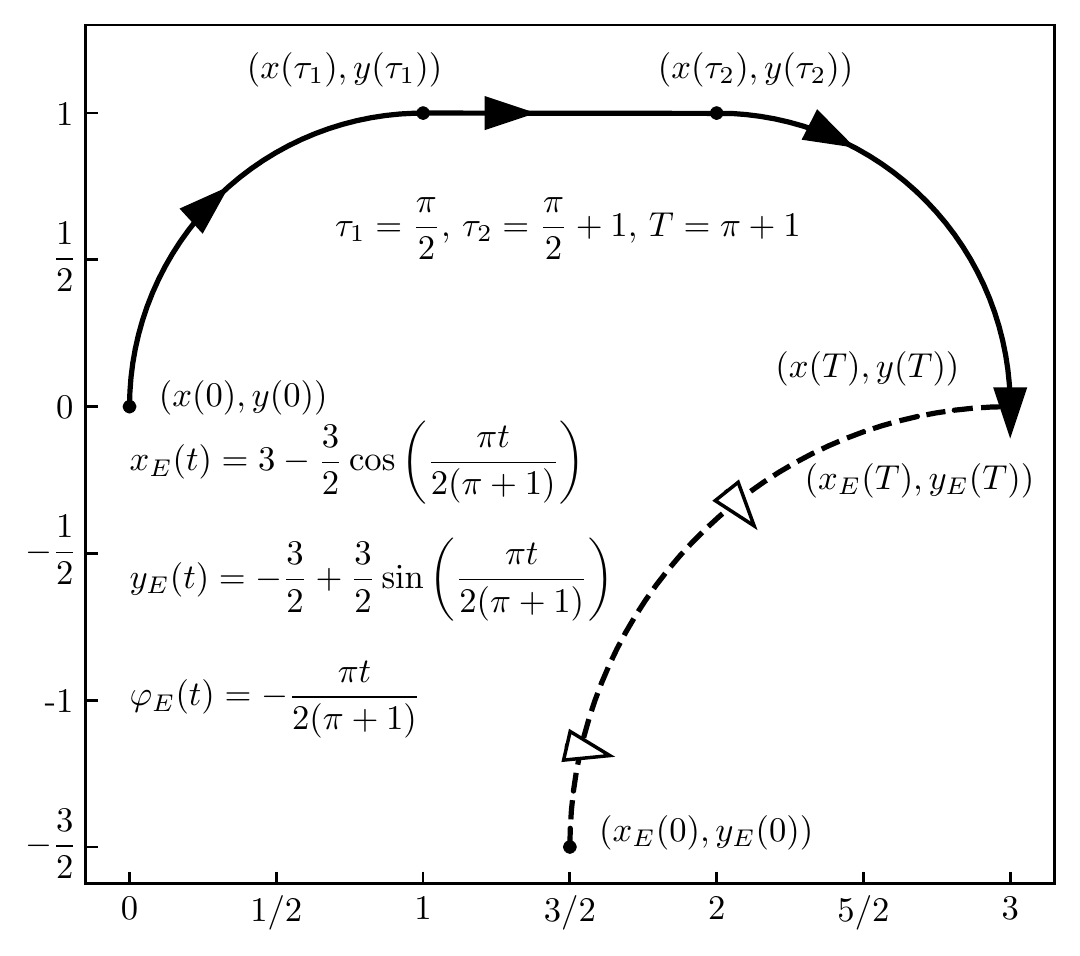}
        \caption{RSR-optimal interception case ($s = -1$, $\sigma = -1$).}
        \label{fig:rSr}
    \end{center}
\end{figure}
\begin{figure}[h!]
    \begin{center}
        \includegraphics[width=\linewidth]{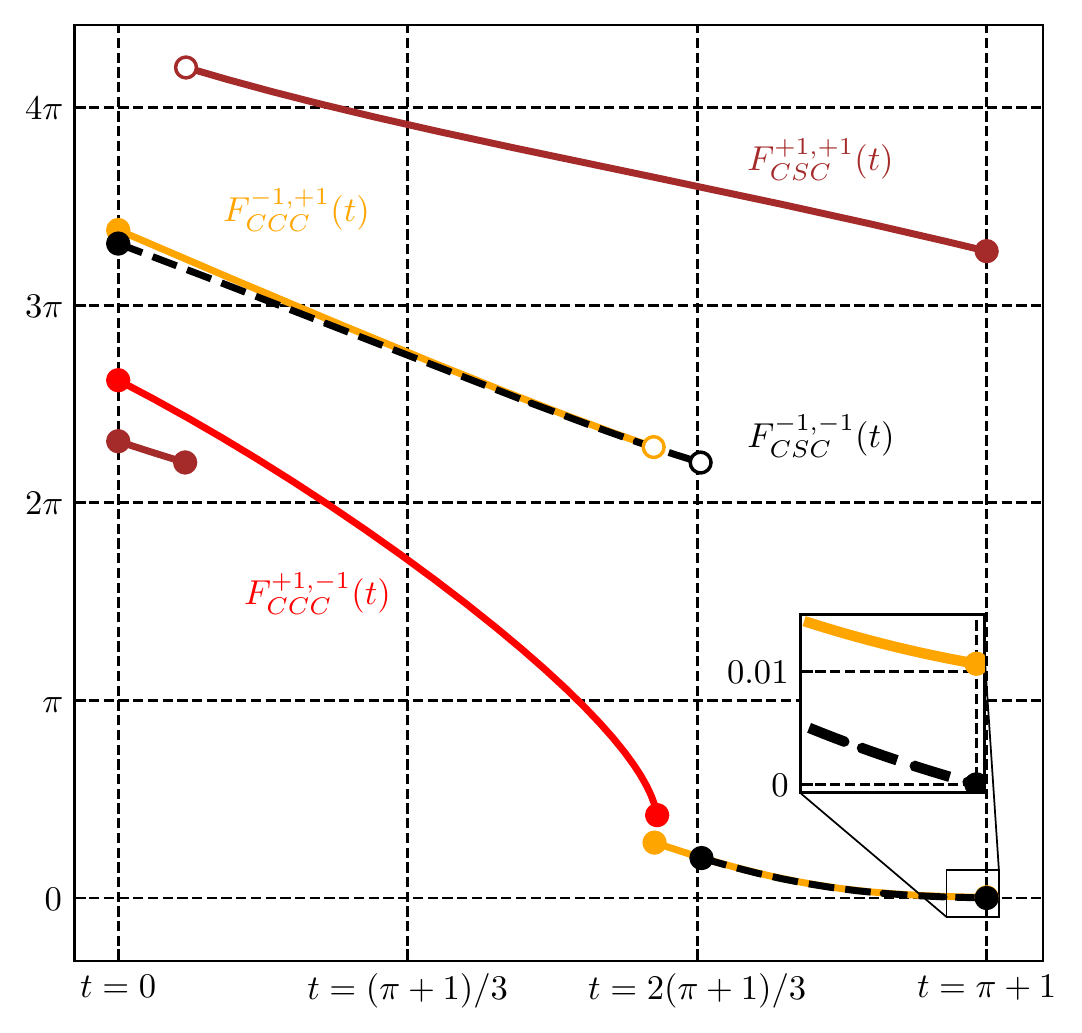}
        \caption{Graphs of functions $F^{+1, +1}_{CSC}$, $F^{-1, -1}_{CSC}$, $F^{-1, +1}_{CCC}$, $F^{+1, -1}_{CCC}$: RSR-optimal interception case.}
        \label{fig:rSr_F1}
    \end{center}
\end{figure}
\begin{figure}[h!]
    \begin{center}
        \includegraphics[width=\linewidth]{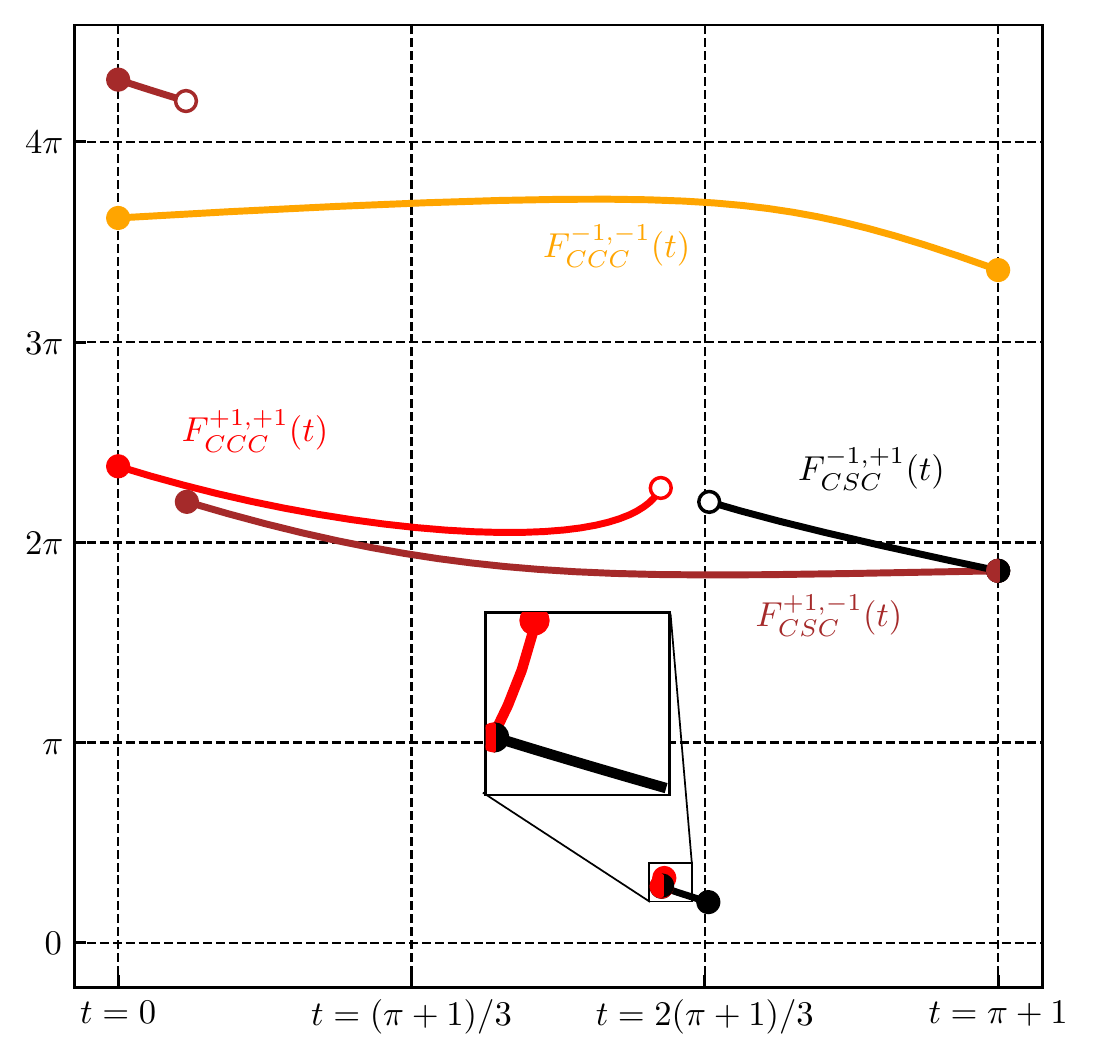}
        \caption{Graphs of functions $F^{+1, -1}_{CSC}$, $F^{-1, +1}_{CSC}$, $F^{+1, +1}_{CCC}$, $F^{-1, -1}_{CCC}$: RSR-optimal interception case.}
        \label{fig:rSr_F2}
    \end{center}
\end{figure}

Consider the case of the circular movement of the target as shown in Fig.~\ref{fig:rSr} by a dashed line. The analysis of the equations \eqref{eq:optimal_T_for_CSC}, \eqref{eq:optimal_T_for_CCC} is presented in Figs.~\ref{fig:rSr_F1},~\ref{fig:rSr_F2}. In these figures one can see that $t = \pi + 1$ is a minimal time moment when the function $F^{-1, -1}_{CSC}$ becomes equal to zero. The other function graphs in the figures do not cross the horizontal axis zero level on the section \mbox{$[0, \pi + 1]$}. Since $T = \pi + 1 < 2\pi$ we should not analyse the cases of the presence of a cycle. Hence, the RSR-trajectory ($s = -1$, $\sigma = -1$) is an optimal trajectory. The mirror symmetry allows to state that the LSL-trajectory ($s = +1$, $\sigma = +1$) is able to be optimal for the $E^*$ target trajectory.

\begin{figure}[h!]
    \begin{center}
        \includegraphics[width=\linewidth]{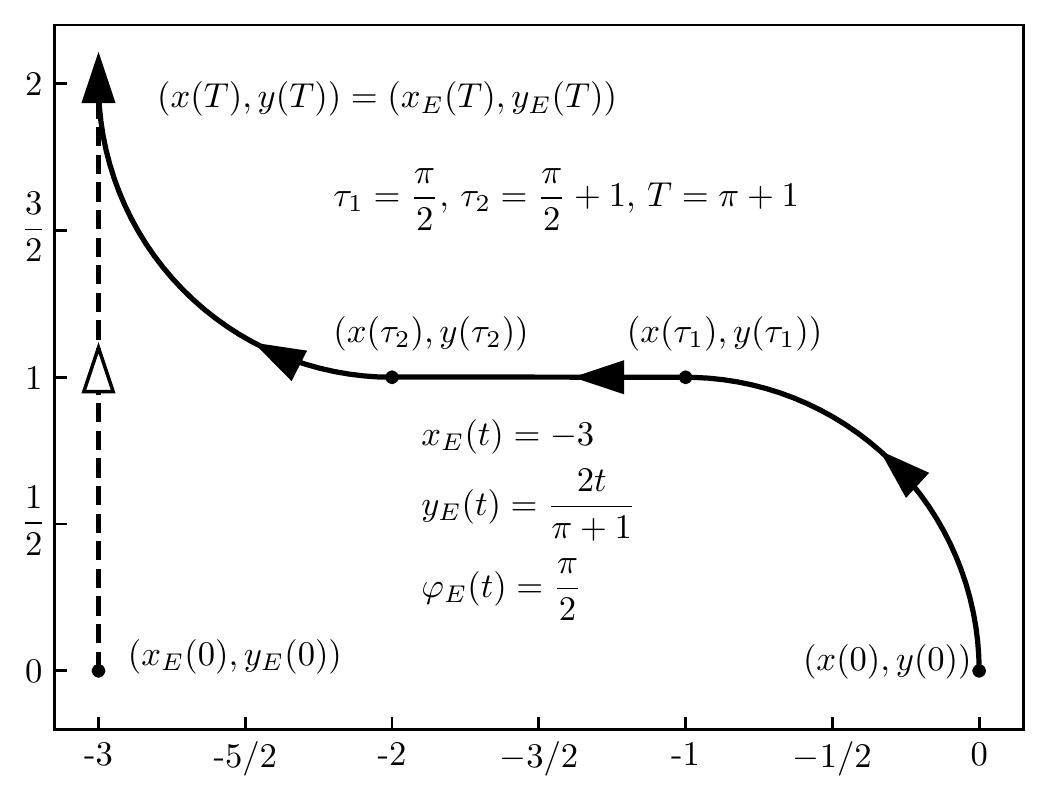}
        \caption{LSR-optimal interception case ($s = +1$, $\sigma = -1$).}
        \label{fig:lSr}
    \end{center}
\end{figure}

\begin{figure}[h!]
    \begin{center}
        \includegraphics[width=\linewidth]{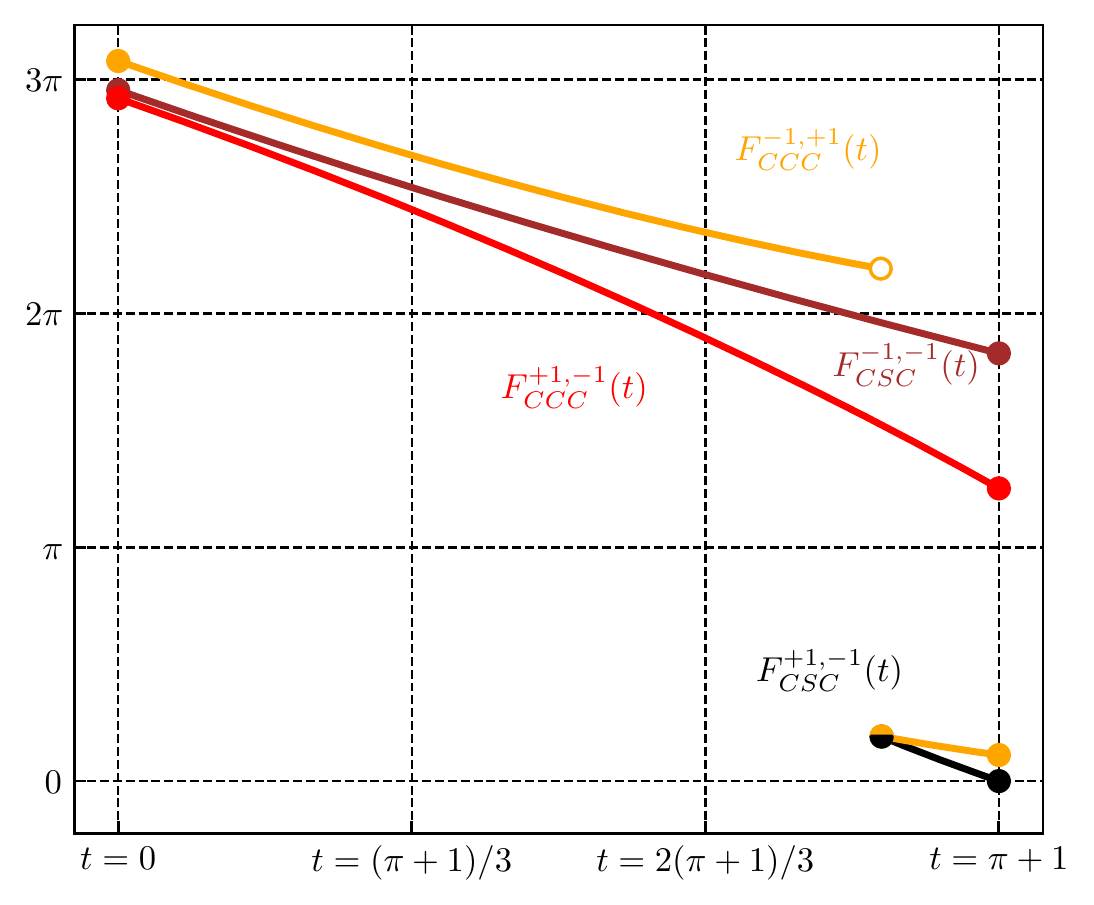}
        \caption{Graphs of functions $F^{-1, +1}_{CCC}$, $F^{-1, -1}_{CSC}$, $F^{+1, -1}_{CCC}$, $F^{+1, -1}_{CSC}$: LSR-optimal interception case.}
        \label{fig:lSr_F1}
    \end{center}
\end{figure}

\begin{figure}[h!]
    \begin{center}
        \includegraphics[width=\linewidth]{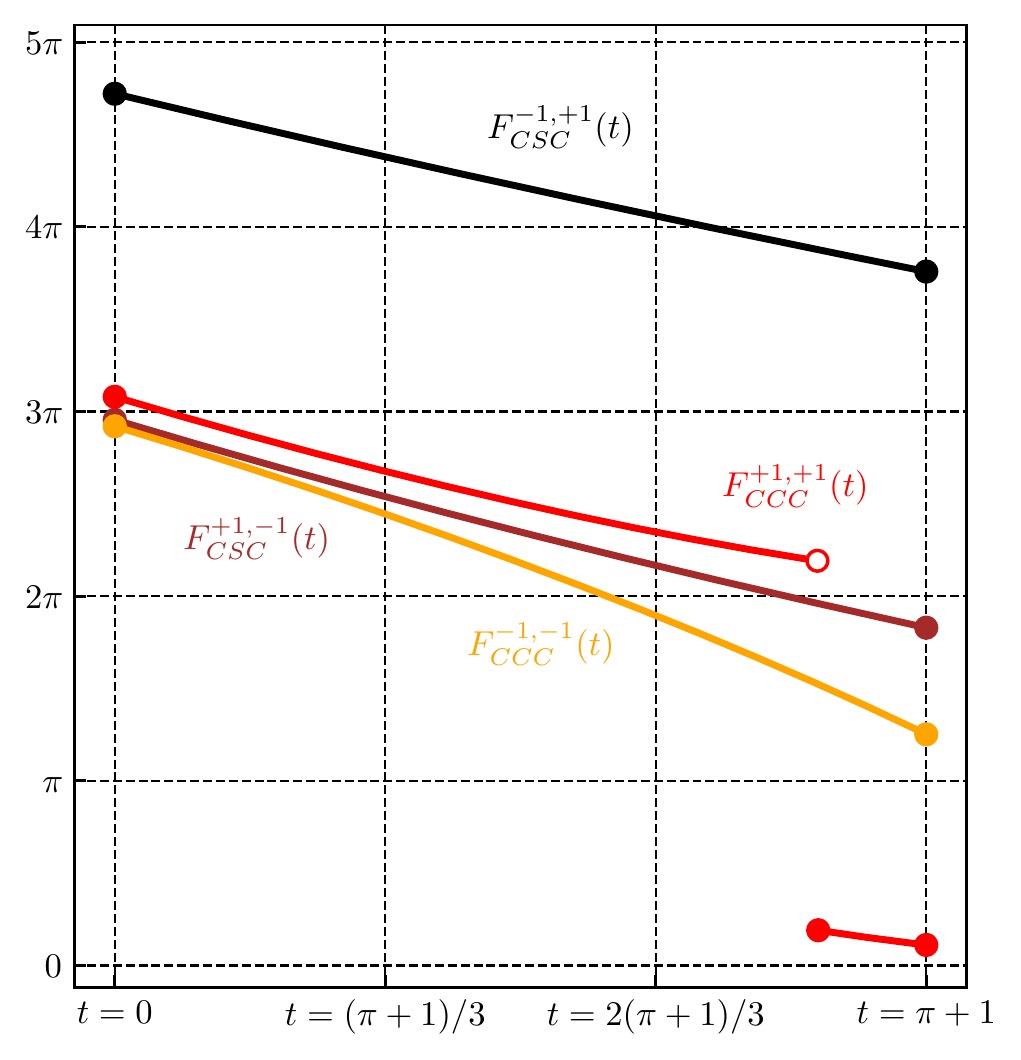}
        \caption{Graphs of functions $F^{-1, +1}_{CSC}$, $F^{+1, +1}_{CCC}$, $F^{+1, -1}_{CSC}$, $F^{-1, -1}_{CCC}$: LSR-optimal interception case.}
        \label{fig:lSr_F2}
    \end{center}
\end{figure}

Next, consider the case of the uniform rectilinear movement of the target (see Fig.~\ref{fig:lSr}). The analysis of the equations \eqref{eq:optimal_T_for_CSC}, \eqref{eq:optimal_T_for_CCC} in Figs.~\ref{fig:lSr_F1},~\ref{fig:lSr_F2} demonstrates that the function $F^{+1, -1}_{CSC}$ becomes equal to zero at $t = \pi + 1$ and the other functions in the figures do not. Similarly, since $T = \pi + 1 < 2\pi$ we should not analyse the cases of the presence of a cycle. Hence, the LSR-trajectory ($s = +1$, $\sigma = -1$) is an optimal trajectory. The mirror symmetry allows to state that the RSL-trajectory ($s = -1$, $\sigma = +1$) is able to be optimal for the $E^*$ target trajectory.


\begin{figure}[h!]
    \begin{center}
        \includegraphics[width=\linewidth]{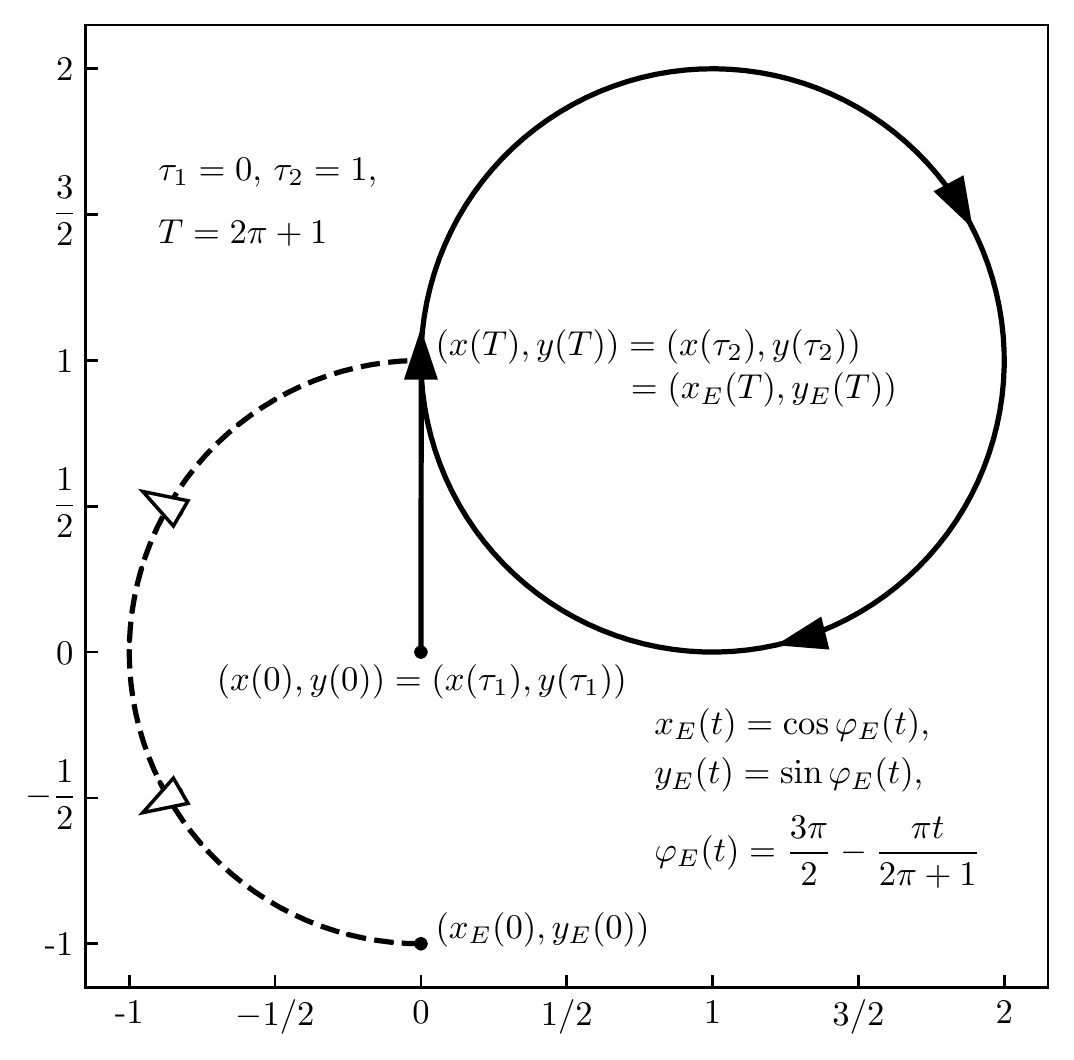}
        \caption{SC-optimal interception case ($\sigma = -1$).}
        \label{fig:SC}
    \end{center}
\end{figure}

\begin{figure}[h!]
    \begin{center}
        \includegraphics[width=\linewidth]{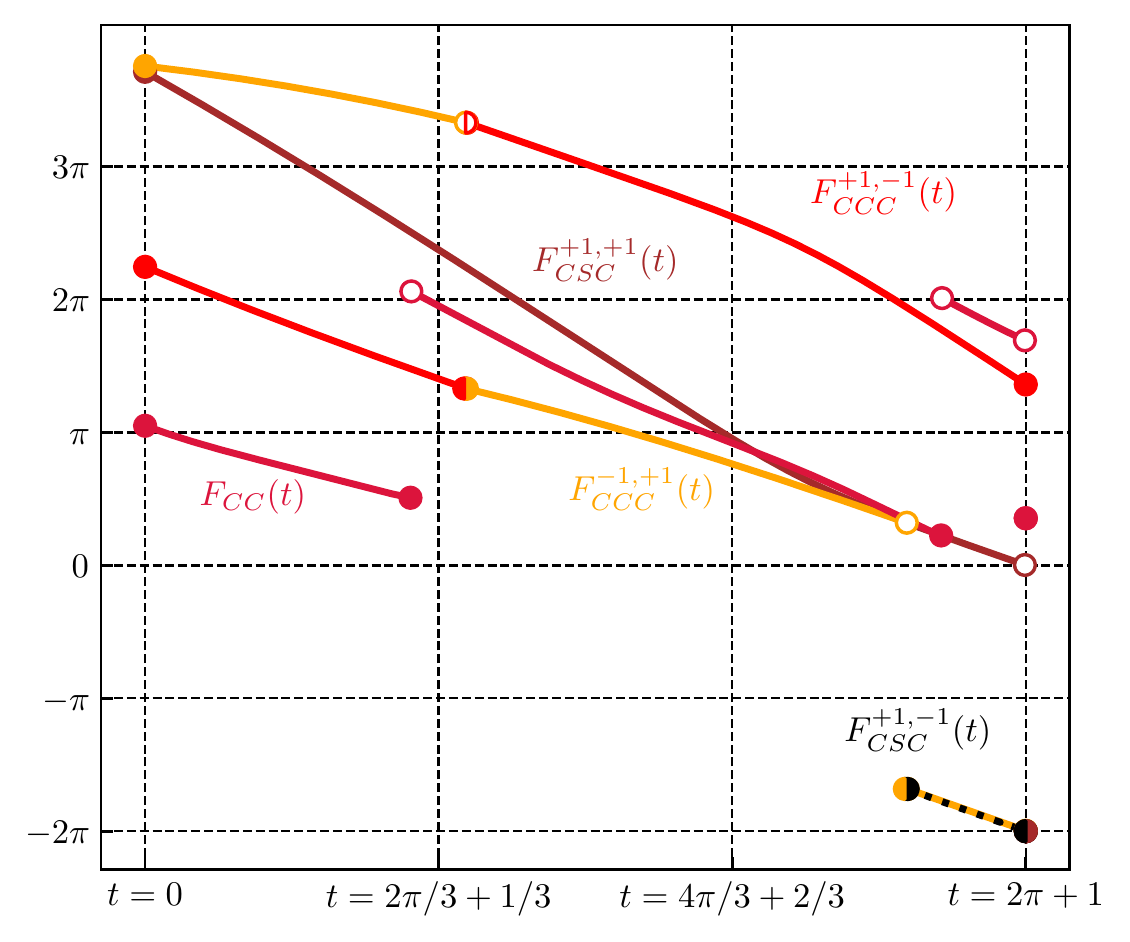}
        \caption{Graphs of functions $F^{+1, -1}_{CCC}$, $F^{+1, +1}_{CSC}$, $F_{CC}$, $F^{-1, +1}_{CCC}$, $F^{+1, -1}_{CSC}$: SC-optimal interception case.}
        \label{fig:SC_F1}
    \end{center}
\end{figure}

\begin{figure}[h!]
    \begin{center}
        \includegraphics[width=\linewidth]{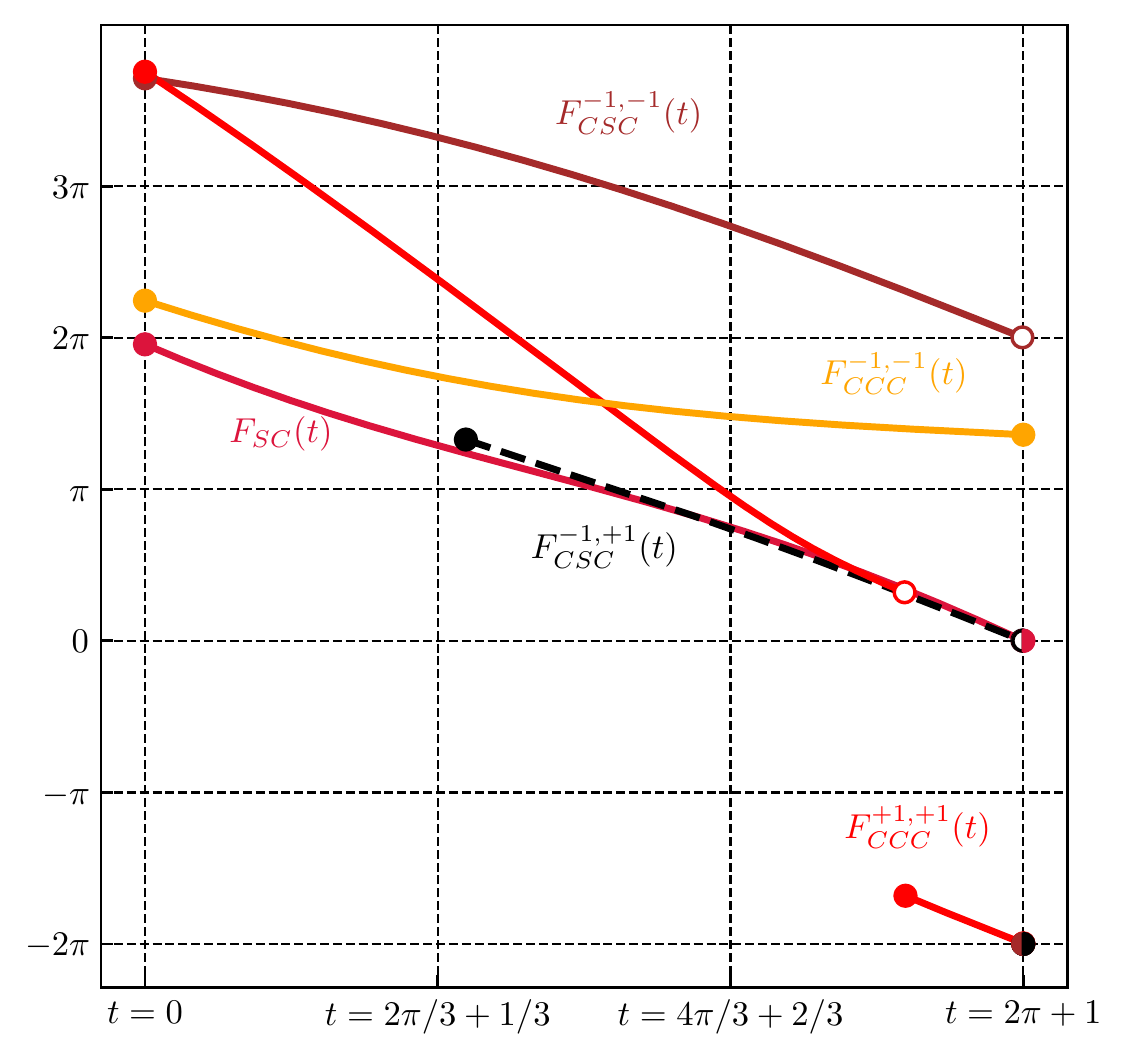}
        \caption{Graphs of functions $F^{-1, -1}_{CSC}$, $F^{-1, -1}_{CCC}$, $F_{SC}$, $F^{-1, +1}_{CSC}$, $F^{+1, +1}_{CCC}$: SC-optimal interception case.}
        \label{fig:SC_F2}
    \end{center}
\end{figure}

Finally, consider the case of the circular movement of the target (see Fig.~\ref{fig:SC}). The analysis of the equations \eqref{eq:optimal_T_for_CSC}, \eqref{eq:optimal_T_for_CSC_with_cycle}, \eqref{eq:optimal_T_for_CCC}, \eqref{eq:optimal_T_for_CCC_with_cycle} in Figs.~\ref{fig:SC_F1},~\ref{fig:SC_F2} reveals that the function $F_{SC}$ becomes equal to zero at $t = 2\pi + 1$ and the other functions do not. In this case all 10 candidates have to be analysed since $T = 2\pi + 1 > 2\pi$. As it is seen in Fig.~\ref{fig:SC}, the optimal interception trajectory of the Dubins car contains a cycle and so as it is mentioned above the parameters $(s,\sigma)$ may be arbitrary chosen. 


\subsection{CCC-optimal examples}

\begin{figure}
    \begin{center}
        \includegraphics[width=0.7\linewidth]{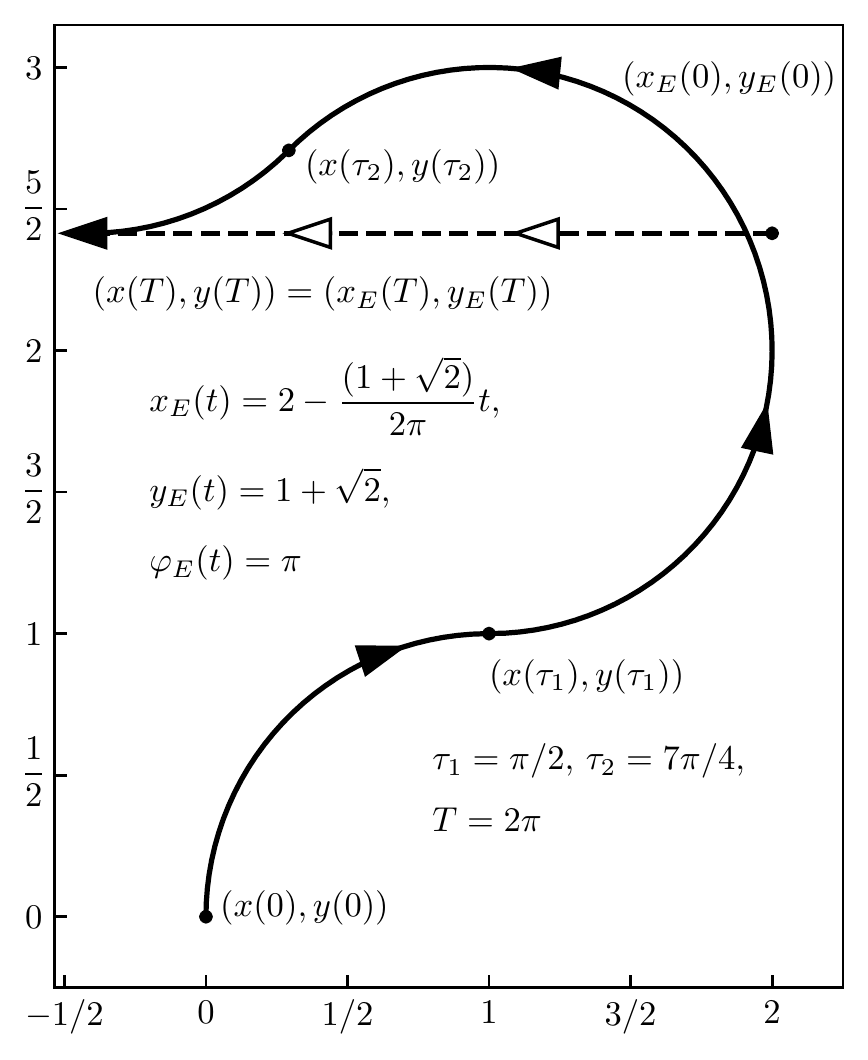}
        \caption{RLR-optimal interception case ($s = -1$, $\mu = -1$).}
        \label{fig:rLr}
    \end{center}
\end{figure}

\begin{figure}
    \begin{center}
        \includegraphics[width=\linewidth]{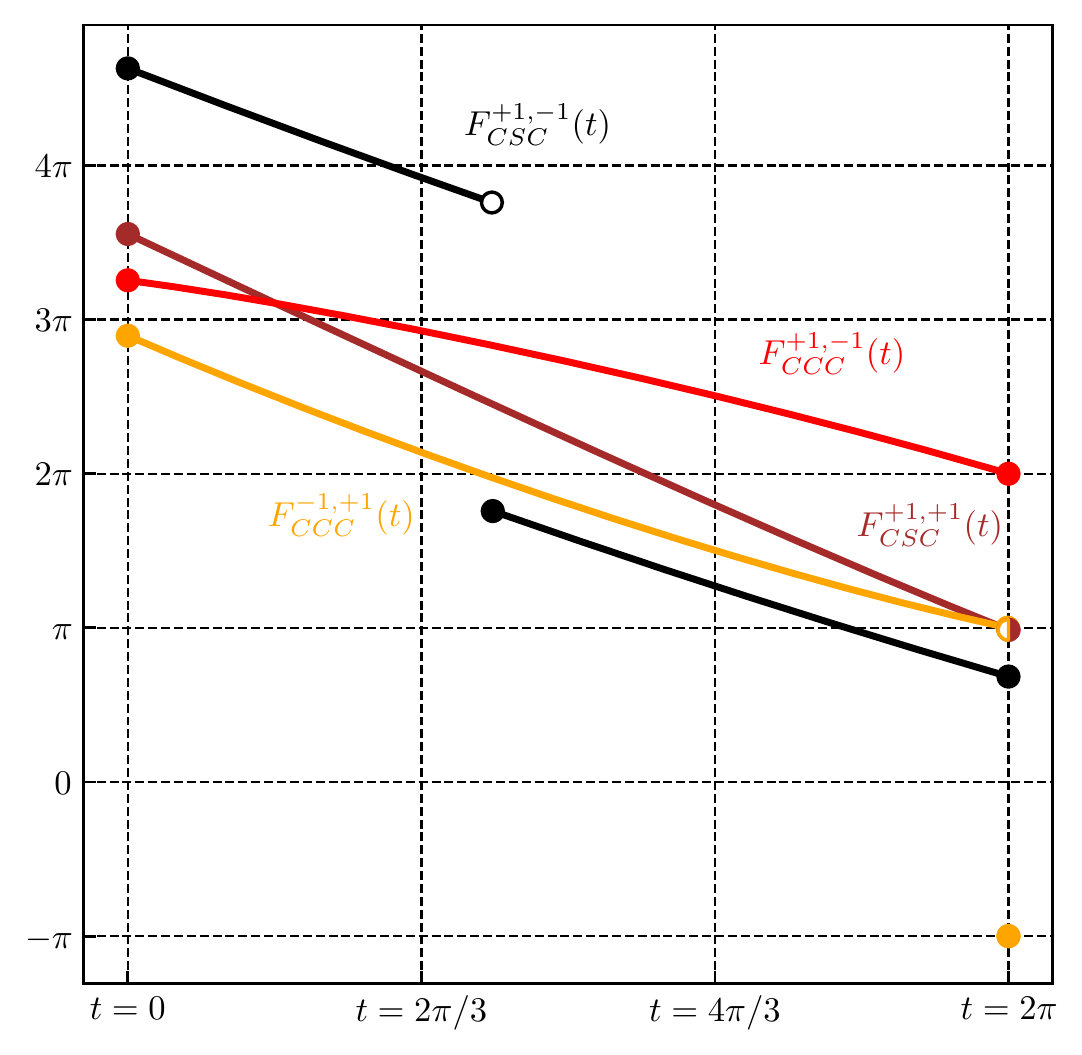}
        \caption{Graphs of functions $F^{+1, -1}_{CSC}$, $F^{+1, -1}_{CCC}$, $F^{-1, +1}_{CCC}$, $F^{+1, +1}_{CSC}$: RLR-optimal interception case ($\mu = -1$).}
        \label{fig:rLr_F1}
    \end{center}
\end{figure}

\begin{figure}
    \begin{center}
        \includegraphics[width=\linewidth]{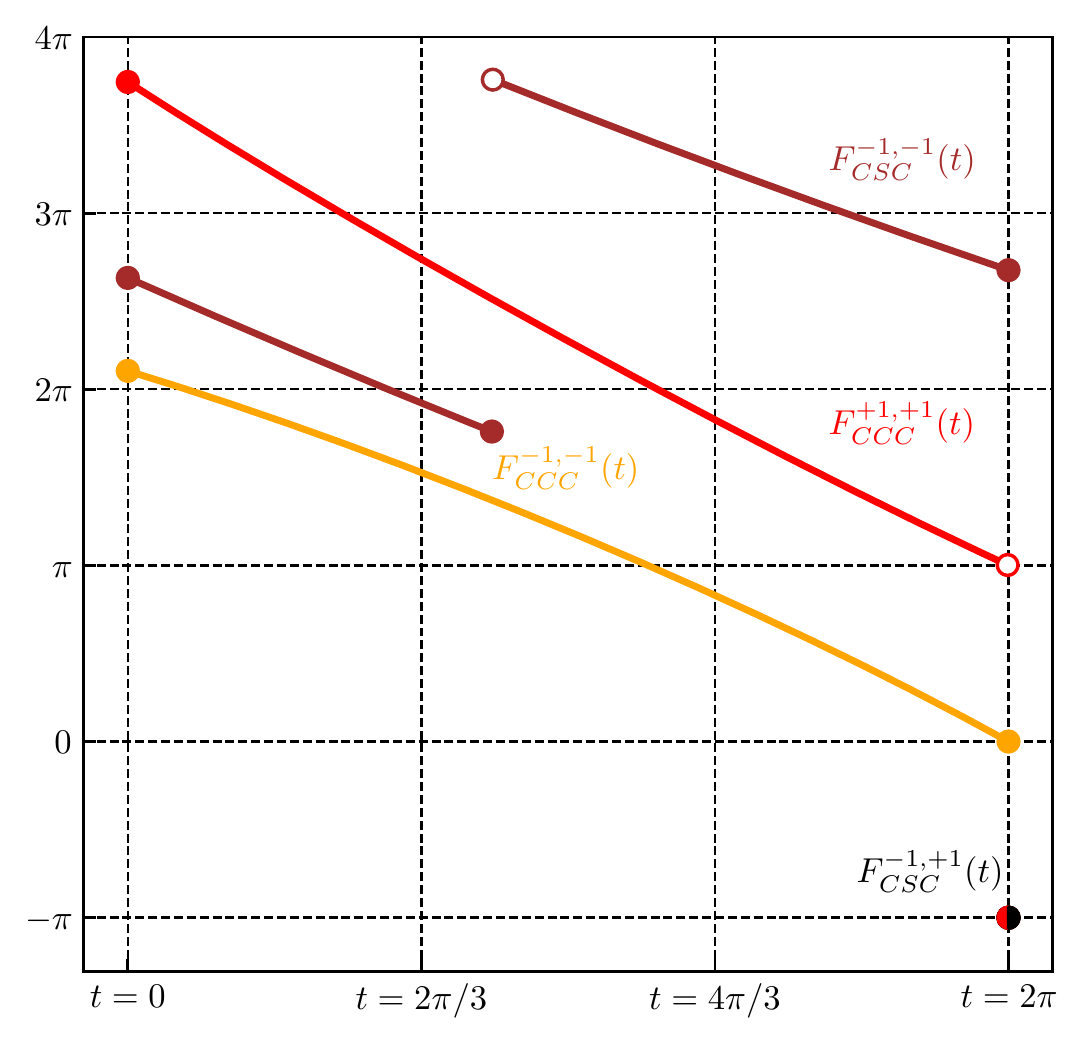}
        \caption{Graphs of functions $F^{-1, -1}_{CSC}$, $F^{+1, +1}_{CCC}$, $F^{-1, -1}_{CCC}$, $F^{-1, +1}_{CSC}$: RLR-optimal interception case ($\mu = -1$).}
        \label{fig:rLr_F2}
    \end{center}
\end{figure}

Consider the other case of the uniform rectilinear movement of the target (see Fig.~\ref{fig:rLr}). The analysis of the equations \eqref{eq:optimal_T_for_CSC}, \eqref{eq:optimal_T_for_CCC} in Figs.~\ref{fig:rLr_F1},~\ref{fig:rLr_F2} demonstrates that the function $F^{-1, -1}_{CCC}$ becomes equal to zero at time $t = 2\pi$ and the other functions in the figures do not. Since the target is not located at $(0, 0, \pi/2)$ at time $t = 2\pi$ we should not analyse the cases of the presence of a cycle. Therefore, the RLR-trajectory ($s = -1$, $\mu = -1$) is an optimal trajectory. The mirror symmetry allows to state that the LRL-trajectory ($s = +1$, $\mu = -1$) is able to be optimal for the $E^*$ target trajectory.


\begin{figure}
    \begin{center}
        \includegraphics[width=0.7\linewidth]{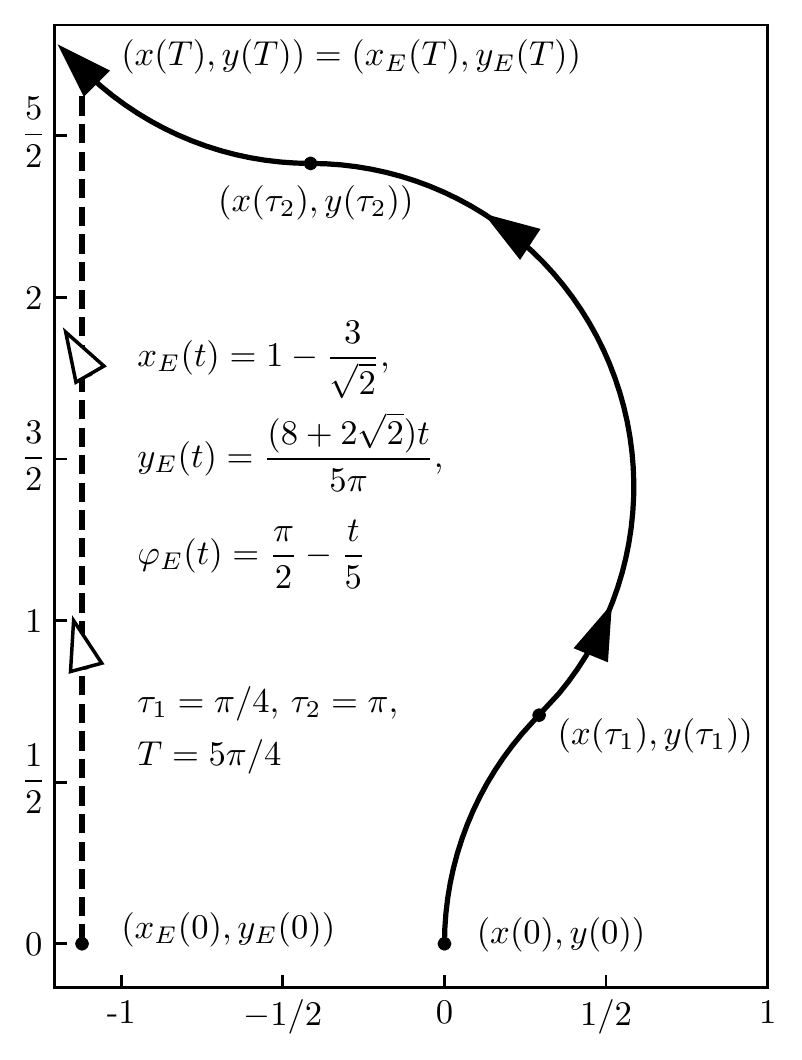}
        \caption{RLR-optimal interception case ($s = -1$, $\mu = +1$).}
        \label{fig:rlr_u}
    \end{center}
\end{figure}

\begin{figure}
    \begin{center}
        \includegraphics[width=\linewidth]{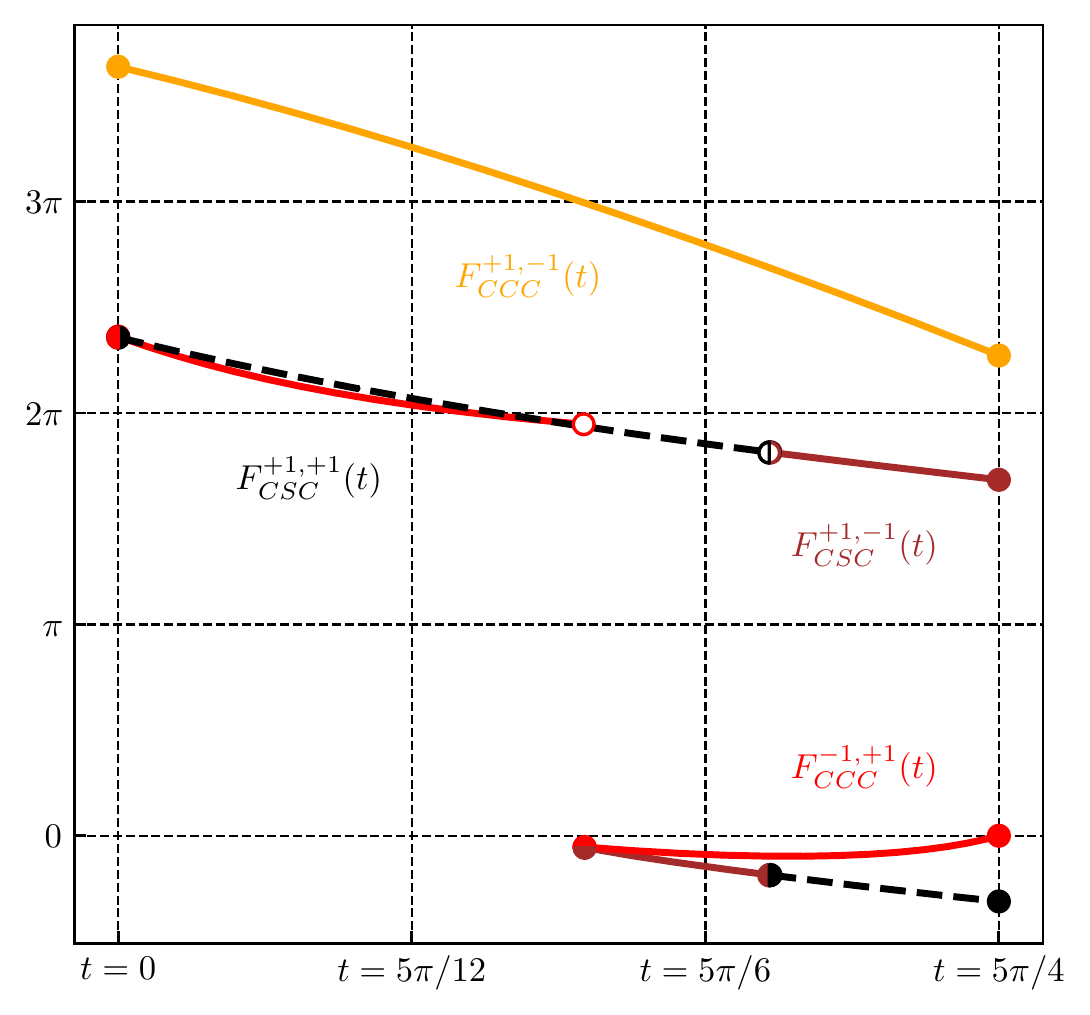}
        \caption{Graphs of functions $F^{+1, -1}_{CCC}$, $F^{+1, +1}_{CSC}$, $F^{+1, -1}_{CSC}$, $F^{-1, +1}_{CCC}$: RLR-optimal interception case ($\mu = +1$).}
        \label{fig:rlr_u_F1}
    \end{center}
\end{figure}

\begin{figure}
    \begin{center}
        \includegraphics[width=\linewidth]{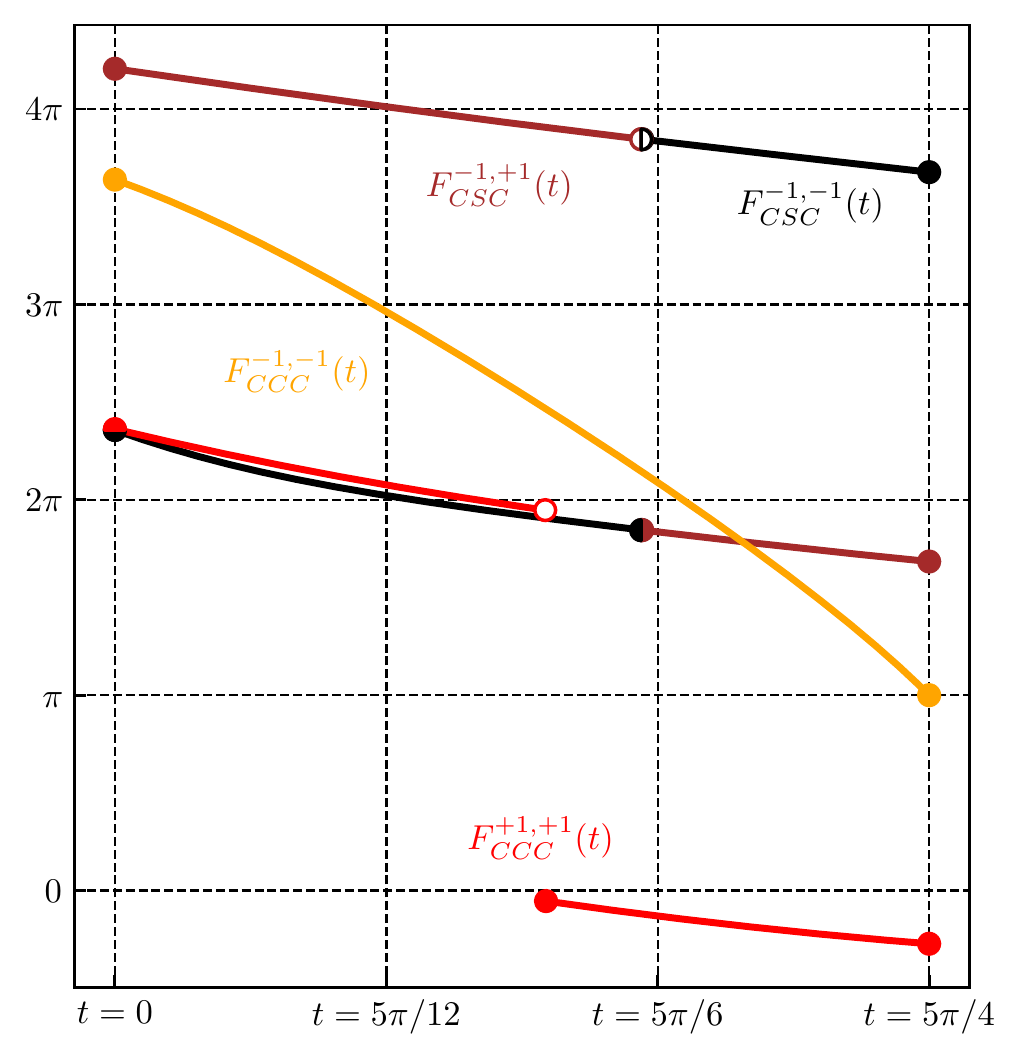}
        \caption{Graphs of functions $F^{-1, +1}_{CSC}$, $F^{-1, -1}_{CSC}$, $F^{-1, -1}_{CCC}$, $F^{+1, +1}_{CCC}$: RLR-optimal interception case ($\mu = +1$).}
        \label{fig:rlr_u_F2}
    \end{center}
\end{figure}

Next, consider the case of the uniform rectilinear movement of the target (see Fig.~\ref{fig:rlr_u}). The target orientation $\varphi_E(t)$ changes clockwise, uniformly and starts from $\varphi_E(0) = \pi/2$. The analysis of the equations \eqref{eq:optimal_T_for_CSC}, \eqref{eq:optimal_T_for_CCC} in Figs.~\ref{fig:rlr_u_F1},~\ref{fig:rlr_u_F2} demonstrates that the function $F^{-1, +1}_{CCC}$ becomes equal to zero at time $t = 5\pi/4$ and the other functions in the figures do not. Since $T = 5\pi/4 < 2\pi$ we may analyse only the cycle-free cases. Hence, the RLR-trajectory ($s = -1$, $\mu = +1$) is an optimal trajectory. The mirror symmetry allows to state that the LRL-trajectory ($s = +1$, $\mu = +1$) is able to be optimal for the $E^*$ target trajectory.


\begin{figure}
    \begin{center}
        \includegraphics[width=\linewidth]{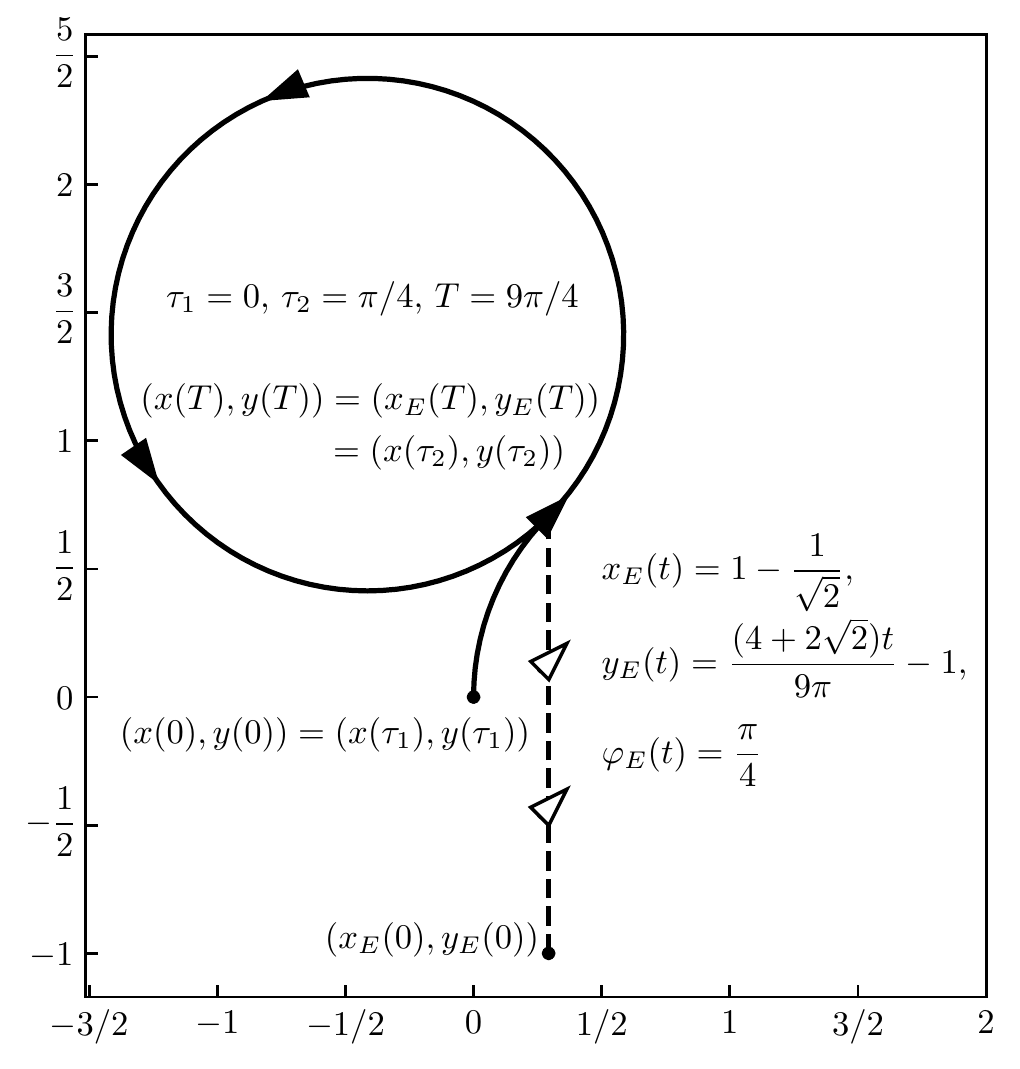}
        \caption{CC-optimal interception case ($s = +1$).}
        \label{fig:Cc}
    \end{center}
\end{figure}

\begin{figure}
    \begin{center}
        \includegraphics[width=\linewidth]{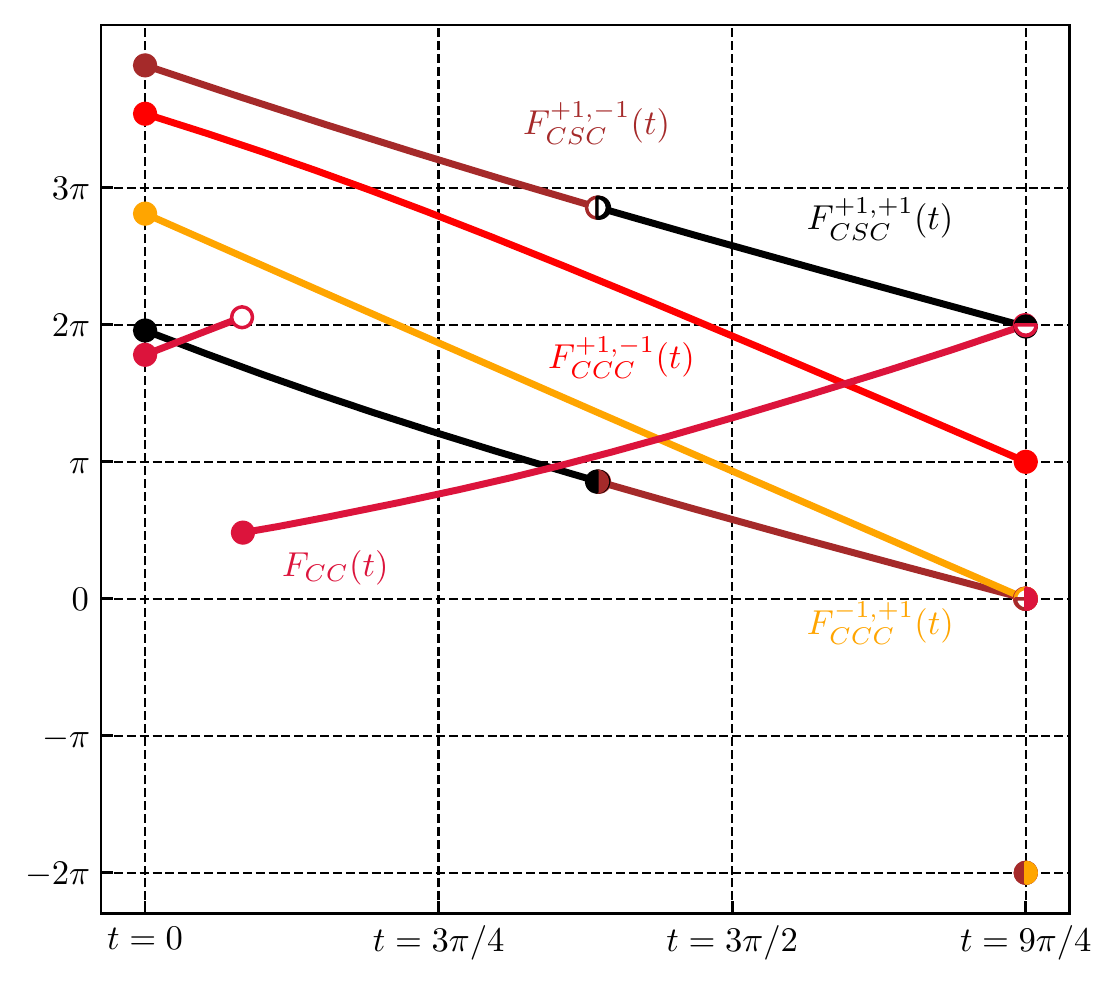}
        \caption{Graphs of functions $F^{+1, -1}_{CSC}$, $F^{+1, +1}_{CSC}$, $F^{+1, -1}_{CCC}$, $F_{CC}$, $F^{-1, +1}_{CCC}$: CC-optimal interception case.}
        \label{fig:Cc_F1}
    \end{center}
\end{figure}

\begin{figure}
    \begin{center}
        \includegraphics[width=\linewidth]{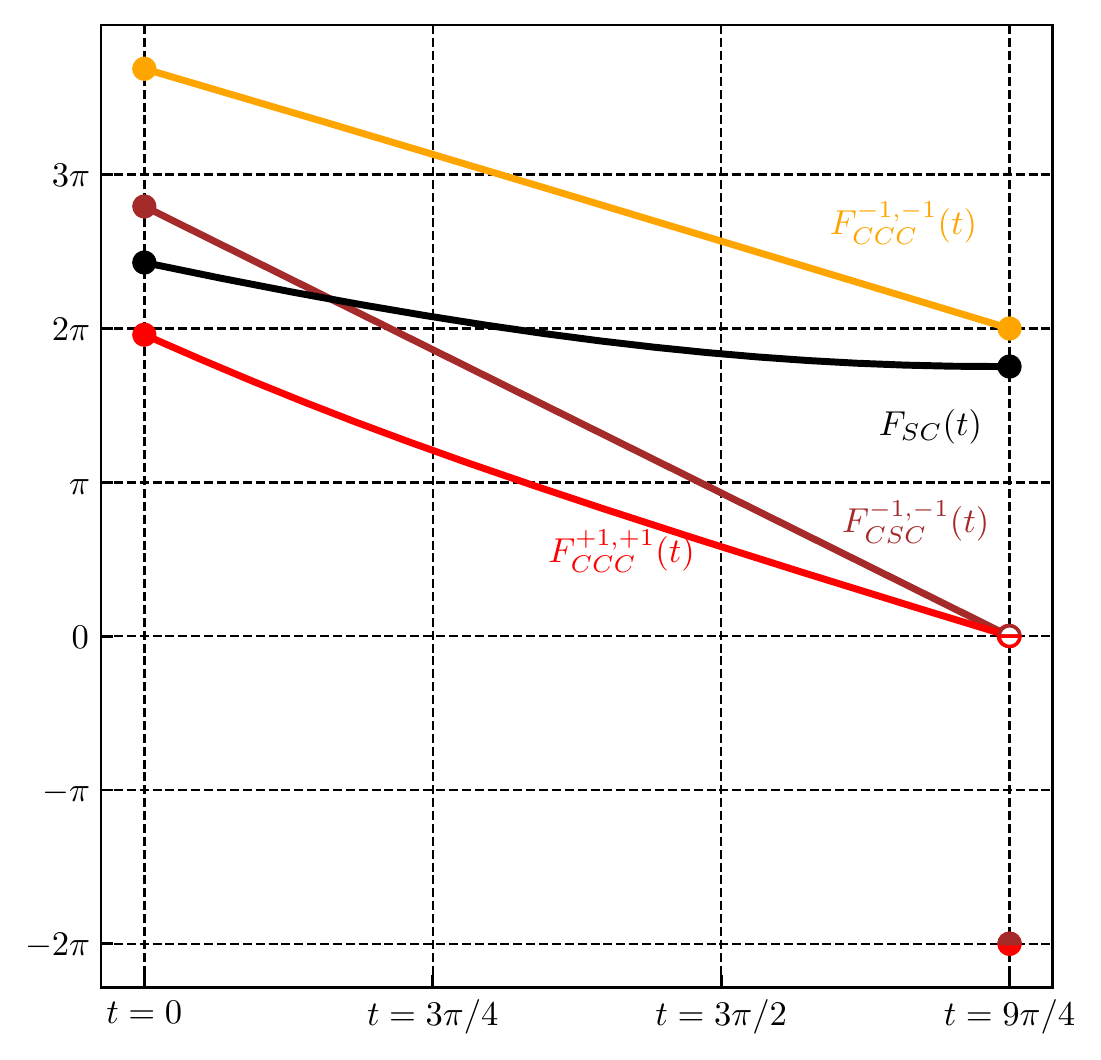}
        \caption{Graphs of functions $F^{-1, -1}_{CCC}$, $F_{SC}$, $F^{-1, -1}_{CSC}$, $F^{+1, +1}_{CCC}$: CC-optimal interception case.}
        \label{fig:Cc_F2}
    \end{center}
\end{figure}

Finally, consider the case of the uniform rectilinear movement of the target that is shown in Fig.~\ref{fig:Cc}. The analysis of the equations \eqref{eq:optimal_T_for_CSC}, \eqref{eq:optimal_T_for_CSC_with_cycle}, \eqref{eq:optimal_T_for_CCC}, \eqref{eq:optimal_T_for_CCC_with_cycle} in Figs.~\ref{fig:Cc_F1},~\ref{fig:Cc_F2} demonstrates that the function $F_{CC}$ becomes equal to zero at time $t = 9\pi/4$ and the other functions in the figures do not.


\section{Conclusions}\label{sec:conclusion}

For the problem of the time-minimal lateral interception of a predefined moving target by the Dubins car we obtained algebraic equations, the solution of which allows one to calculate all the parameters of the optimal control. Since the target movement is specified in the form of an arbitrary continuous vector-function $E = (x_E, y_E, \varphi_E)$, the obtained analytical results are suitable for a wide range of cases of the target movement. Moreover, the obtained results are applicable to the problem of the time-minimal achievement of the desired state for the Dubins car moving  under time-depended wind-conditions. 

Instead of the maximum principle application, we used another methodological approach which deals with geometric properties of the reachable set $\mathcal{R}(t)$ and the corresponding multi-valued mapping. This technique can be useful for considering interception problems by mobile vehicles whose dynamics differ from the dynamics of the Dubins car. Theorem~\ref{th:optimal_interception_point} determines the optimal interception location on the reachable set, and the analytical description of the multi-valued mapping $\mathcal{E}$ makes it possible to reduce the initial problem to solving algebraic equations. Theorem~\ref{th:CSC_CCC_additional_restriction} gives the additional constraints on the parameters of the optimal control without losing the optimality of the solution and reduces the number of algebraic equations required for consideration to 10. The described examples prove the necessity to consider all 10 equations in general case.

Topics to be addressed in future work include an investigation of effective numerical methods for solving the obtained algebraic transcendental equations. These algebraic equations are discontinuous and not defined on the entire number line, which raises a sophisticated challenge for the standard numerical methods.

\bibliographystyle{plain}
\bibliography{root}

\end{document}